\documentclass[12pt, twoside, reqno]{amsart}
\usepackage{amssymb,amscd,amsxtra,latexsym,bbm,graphicx,epsfig,epic,
eepic,oldgerm,bm,psfrag,amsthm,mathrsfs,bm,bbm}
\usepackage{caption}
\usepackage{subcaption}
\usepackage{tikz-cd}
\usepackage[pdftex]{hyperref}

\input xy
\xyoption{all}
\CompileMatrices
    
\theoremstyle{plain}
\newtheorem{theorem}{Theorem}[section]
\newtheorem*{theorem*}{Theorem}
\newtheorem{lemma}[theorem]{Lemma}
\newtheorem{proposition}[theorem]{Proposition}

\newtheorem*{remark*}{Remark}
\newtheorem*{remarks*}{Remarks}
\newtheorem{remark}[theorem]{Remark}

\newtheorem{example}[theorem]{Example}
\newtheorem*{example*}{Example}
\newtheorem*{examples*}{Examples}
\newtheorem{definition}[theorem]{Definition}
\newtheorem*{definition*}{Definition}

\newtheorem{question}[theorem]{Question}

\newtheorem{assumption}[theorem]{Assumption}

\newcommand{\proofend}{\hspace*{\fill} $\Box$\\}

\newcommand{\ign}[1]{}

\def\1{\:\!}
\def\2{\;\!}

\def\Diff{\operatorname{Diff}}

\def\Diffc0{\operatorname{Diff^c_0}}

\def\Sympc0{\operatorname{Symp^c_0}}

\def\idd{\operatorname{id}}

\def\GL{\operatorname{GL}}

\def\Fix{\operatorname{Fix}}
\def\maslov{\operatorname{Maslov}}

\def\momega{(M,\omega)}
\def\momegam{(M,\omega,\mu)}

\def\CC{\mathbbm{C}}

\def\RR{\mathbbm{R}}

\def\ZZ{\mathbbm{Z}}

\def\RP{\operatorname{\mathbbm{R}P}}
\def\CP{\operatorname{\mathbbm{C}P}}

\def\Diff{\operatorname{Diff}}

\def\DP{\operatorname{DP}}

\def\ni{\noindent}

\def\.{\mskip1mu}
\def\?{\mskip-1mu}

\renewcommand{\colon}{\nobreak\mskip2mu\mathpunct{}\nonscript\mkern-\thinmuskip{:}\mskip6muplus1mu\relax}

\def\proof{\noindent {\it Proof. \;}}
\newcommand{\proofof}[1]{\ni {\it Proof of #1. }}

\begin{document}

\title[]{Real Lagrangian tori and versal deformations}

\author{Jo\'e Brendel}  
\thanks{partially supported by SNF grant 200020-144432/1.}
\address{Jo\'e Brendel,
Institut de Math\'ematiques,
Universit\'e de Neuch\^atel}
\email{joe.brendel@unine.ch}


\date{\today}
\thanks{2010 {\it Mathematics Subject Classification.}
Primary~, Secondary~}

\maketitle

\begin{abstract}
Can a given Lagrangian submanifold be realized as the fixed point set of an anti-symplectic involution? If so, it is called \emph{real}. We give an obstruction for a closed Lagrangian submanifold to be real in terms of the displacement energy of nearby Lagrangians. Applying this obstruction to toric fibres, we obtain that the central fibre of many (and probably all) toric monotone symplectic manifolds is real only if the corresponding moment polytope is centrally symmetric. Furthermore, we embed the Chekanov torus in all toric monotone symplectic manifolds and show that it is exotic and not real, extending Kim's result \cite{Kim19b} for $S^2 \times S^2$. Inside products of $S^2$, we show that all products of Chekanov tori are pairwise distinct and not real either. These results indicate that real tori are rare.

Our methods are elementary in the sense that we do not use~$J$-holomorphic curves. Instead, we rely on symplectic reduction and the displacement energy of product tori in~$\RR^{2n}$.
\end{abstract}

\section{Introduction}
\label{sec:intro}

A Lagrangian submanifold $L$ in a symplectic manifold $(M,\omega)$ is said to be {\it real} if there is an anti-symplectic involution~$\sigma$ of~$M$ such that $L$ is the fixed point set of~$\sigma$ or a connected component thereof. Here, an involution is a map satisfying $\sigma \circ \sigma = {\rm id}$, and anti-symplectic means that $\sigma^* \omega = - \omega$. An example is the equator of the 2-sphere with its Euclidean area form, which is the fixed point set of the reflection about the equatorial plane, and taking products of this example we get as real Lagrangian the so-called Clifford torus in $\times_n S^2$. For more examples, see Section \ref{sec:versal}.

\emph{Real} or \emph{not real} are symplectic invariants in the following sense: If~$\varphi$ is a symplectomorphism of $(M,\omega)$ and $L$ is the fixed point set of the anti-symplectic involution~$\sigma$, then $\varphi (L)$ is the fixed point set of the anti-symplectic involution $\varphi \circ \sigma \circ \varphi^{-1}$. There are many other reasons to study real Lagrangian submanifolds, some of which we give at the end of this introduction.

In this paper we address the question if a given closed Lagrangian submanifold of a symplectic manifold
is real. 
An obstruction to being real has been given by J. Kim in~\cite{Kim19b}:
If $L$ is real, then the number of $J$-holomorphic discs $u \colon (D^2, \partial D^2) \to (M,L)$ 
of Maslov index~$2$ passing through a generic point in $L$ must be even.
 In this paper we use a different symplectic invariant as obstruction to being real, namely the displacement energy of nearby Lagrangian submanifolds, a tool invented by Chekanov in~\cite{Che96}. While the Lagrangian submanifolds $L$ that we are interested in usually have infinite displacement energy, nearby Lagrangians can be displaced. This leads to the so-called displacement energy germ $S_L \colon (H_1 (L,\RR),0) \to \RR \cup \{\infty\}$. In our basic result, $(M,\omega)$ is any, not necessarily compact, symplectic manifold.

\begin{theorem} \label{t:dis}
Assume that $L$ is a compact real Lagrangian submanifold of $(M,\omega)$. 
Then the displacement energy germ $S_L \colon (H_1(L, \RR),0) \to \RR \cup \{\infty\}$ is even,
$$S_L(-p) = S_L(p).$$
\end{theorem}

In general, it is hard to compute the displacement energy germ of a Lagrangian~$L$. However, for the special class of fibers of toric symplectic manifolds we show in Section~\ref{sec:displacement} that the displacement energy is intimately related to the moment polytope $\Delta$.

\medskip
\paragraph{\bf Application I: Toric fibres.}
Let $(M,\omega)$ be a toric symplectic manifold with moment map $\mu$ and moment polytope $\Delta = \mu(M)$. For all~$x\in \mathring{\Delta}$, the toric fibre $T_x = \mu^{-1}(x)$ is Lagrangian. These Lagrangian tori are especially well-suited to our methods, since they come with a natural versal deformation defined by varying the base point $a \mapsto T_{x+a}$. Hence, we are led to the question of what the displacement energy of toric fibres looks like as a function of the base point. In other words, we want to understand the function
\begin{equation*}
	e_{\Delta} \colon \Delta \rightarrow  \RR \cup \{\infty\}, \quad x \mapsto  e_M(T_x),
\end{equation*}
\noindent 
where $e_M$ denotes displacement energy. If $T_x$ is real, we get by Theorem \ref{t:dis} that the function $e_{\Delta}$ is invariant under central symmetry in a neighbourhood of $x$.

\begin{figure}
	\begin{subfigure}{0.4\textwidth}
		\centering
		\begin{tikzpicture}[scale=0.465]
			\draw [thick,fill=black!15] (-2,-2)--(-2,4)--(4,-2)--(-2,-2);
			\draw[step=2.0,black!30,thin] (-2.8,-2.8) grid (4.8,4.8);
			\draw[thick] (-2,-2)--(-2,4)--(4,-2)--(-2,-2);
			\draw[thin] (-1.75,-1.75)--(-1.75,3.5)--(3.5,-1.75)--(-1.75,-1.75);
			\draw[thin] (-1.5,-1.5)--(-1.5,3)--(3,-1.5)--(-1.5,-1.5);
			\draw[thin] (-1.25,-1.25)--(-1.25,2.5)--(2.5,-1.25)--(-1.25,-1.25);
			\draw[thin] (-1,-1)--(-1,2)--(2,-1)--(-1,-1);
			\draw[thin] (-0.75,-0.75)--(-0.75,1.5)--(1.5,-0.75)--(-0.75,-0.75);
			\draw[thin] (-0.5,-0.5)--(-0.5,1)--(1,-0.5)--(-0.5,-0.5);
			\draw[thin] (-0.25,-0.25)--(-0.25,0.5)--(0.5,-0.25)--(-0.25,-0.25);
			\draw[thick] (-0.15,0)--(0.15,0)--(0.15,0);
			\draw[thick] (0,0.15)--(0,-0.15)--(0,0.15);
		\end{tikzpicture}
	\caption{$M = \CP^2$}
	\end{subfigure}
	\begin{subfigure}{0.4\textwidth}
  		\centering
		\begin{tikzpicture}[scale=0.7]
			\draw [thick,fill=black!15] (-2,-2)--(2,-2)--(2,2)--(-2,2)--(-2,-2);
			\draw[step=2.0,black!30,thin] (-2.5,-2.5) grid (2.5,2.5);
			\draw[thick] (-2,-2)--(2,-2)--(2,2)--(-2,2)--(-2,-2);
			\draw[thin] (-1.75,-1.75)--(1.75,-1.75)--(1.75,1.75)--(-1.75,1.75)--(-1.75,-1.75);
			\draw[thin] (-1.5,-1.5)--(1.5,-1.5)--(1.5,1.5)--(-1.5,1.5)--(-1.5,-1.5);
			\draw[thin] (-1.25,-1.25)--(1.25,-1.25)--(1.25,1.25)--(-1.25,1.25)--(-1.25,-1.25);
			\draw[thin] (-1,-1)--(1,-1)--(1,1)--(-1,1)--(-1,-1);
			\draw[thin] (-0.75,-0.75)--(0.75,-0.75)--(0.75,0.75)--(-0.75,0.75)--(-0.75,-0.75);
			\draw[thin] (-0.5,-0.5)--(0.5,-0.5)--(0.5,0.5)--(-0.5,0.5)--(-0.5,-0.5);
			\draw[thin] (-0.25,-0.25)--(0.25,-0.25)--(0.25,0.25)--(-0.25,0.25)--(-0.25,-0.25);
			\draw[thick] (-0.1,0)--(0.1,0)--(0.1,0);
			\draw[thick] (0,0.1)--(0,-0.1)--(0,0.1);
		\end{tikzpicture}
  	\caption{$M = S^2\times S^2$}
	\end{subfigure}
\caption{Level sets of the function $e_{\Delta}(x) = e_M(T_x)$.}
\label{fig:level2}
\end{figure}
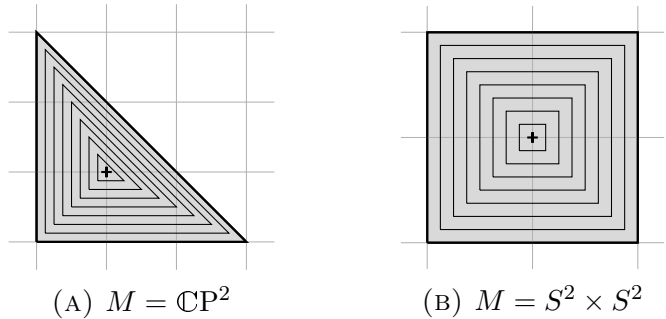

Assume furthermore that $(M,\omega)$ is monotone. In the toric case, this means that we can assume that each facet of the moment polytope lies at affine distance one from the origin, in particular the origin is the only lattice point in the interior. We call the corresponding fibre $T_0$ the \emph{central fibre}. The moment polytope of a toric monotone symplectic manifold is called \emph{monotone}, see \cite{McD11} for details. In this case, the function $e_{\Delta}$ can often be explicitely computed on an open dense subset of $\Delta$ and there is equal to the affine distance to the boundary. In particular the level sets of $e_{\Delta}$ are simply given by rescalings of $\partial\Delta$, see Figure \ref{fig:level2}. As noticed in \cite{BreCheSch20}, this geometric property is implied by the following combinatorial property of the moment polytope: Let $\mathcal{S}(\Delta) = \Delta \cap (-\Delta) \cap \ZZ^n \setminus \{0\}$  be the set of non-zero symmetric lattice points in~$\Delta$. We say that $\Delta$ has property~$FS$ if every facet of~$\Delta$ contains a point of $\mathcal{S}(\Delta)$. This property, which is closely related to the Ewald conjecture, is known to hold for monotone polytopes in dimensions~$n \leq 9$ and is conjectured to hold in all dimensions, in which case requiring property~$FS$ becomes obsolete in all following statements. See Subsection~\ref{ssec:monotone} for a discussion.

Monotonicty has another useful consequence. Since real Lagrangians in monotone symplectic manifolds are automatically monotone as Lagrangian submanifolds, the only candidate to be real among all $T_x$ is the central fibre $T_0$. For this torus we obtain the following.

\begin{theorem} \label{t:toric}
Let $(M,\omega)$ be a toric monotone symplectic manifold whose moment polytope $\Delta$ has property $FS$. If the central fibre $T_0$ is real, then $\Delta$ is centrally symmetric, $\Delta = -\Delta$.
\end{theorem}

Together with J. Kim and J. Moon, we show in \cite{BreKimMoo19} that central symmetry of the moment polytope is a sufficient condition for the central fibre $T_0$ to be real. Under property $FS$, Theorem \ref{t:toric} is therefore an equivalence. For example, the central fibre in $S^2 \times S^2$ is real, whereas the central fibre in $\CP^2$ is not, see Figure \ref{fig:level2}.

\begin{remark}
{\rm
We outline an alternative approach to Theorem~\ref{t:toric} in appendix~\ref{sec:appendix} based on the count of Maslov~$2$ $J$-holomorphic disks with boundary on~$T_0$ which disposes of property~$FS$. This approach was suggested to us by Grigory Mikhalkin and an anonymous referee. 
}
\end{remark}

\medskip
\paragraph{\bf Centrally symmetric polytopes.}
The set of centrally symmetric monotone Delzant polytopes is known. For any natural number $n$ we define the \emph{del Pezzo polytope} $\DP(n) \subset \RR^n$ as the monotone polytope defined by the $2n + 2$ inequalities
\begin{equation*}
	\pm x_1 \leqslant 1 ,\; \pm x_2 \leqslant 1, \ldots , \pm x_n \leqslant 1,\; \pm (x_1 + \ldots + x_n) \leqslant 1.
\end{equation*}
\noindent
For example $\DP(1) = [-1,1]$ and $\DP(2)$ is the moment polytope of the monotone three-fold blow-up of $\CP^2$, see Figure~\ref{fig:2dim} in Section \ref{sec:displacement}. In general, these correspond to two-fold blow-ups of $\times_n S^2$. For $n$ even or $n=1$, the del Pezzo polytopes are centrally symmetric monotone Delzant polytopes. It is thus clear that products of such polytopes~$\DP (n)$ are again centrally symmetric monotone Delzant polytopes. It was proved in~\cite{VosKly84} that the converse is true: The centrally symmetric monotone Delzant polytopes of dimension~$n$ are exactly the products of del Pezzo polytopes $\DP(n_j)$ with $n_j \in \{1,2,4,6,\ldots\}$ and $n = \sum_j n_j$. In order to determine the number $\nu_c(n)$ of centrally symmetric monotone Delzant polytopes in a given dimension $n$, we thus only need to count the number of ways in which $n$ can be written as the sum of ones and even numbers. Let $p(n)$ be the partition function, i.e.\ the function counting the number of ways in which $n$ can be written as the sum of natural numbers. Then for even $n = 2k$,
\begin{equation*}
	\nu_c(2k) = \sum_{j=0}^k p(j).
\end{equation*}
This can be seen as follows. Suppose a decomposition of $2k$ contains~$2m$ ones. Omitting the ones induces a decomposition of $2(k-m)$ into strictly positive even numbers. This is equivalent to a decomposition of $k-m$ into strictly positive integers, whence there are $p(k-m)$ possibilities if the decomposition of $2k$ contains $2m$ ones. Summing over the possible number of ones yields the result. Furthermore $\nu_c(2k+1) = \nu_c(2k)$, since the odd del Pezzo polytopes $\DP(n)$ for $n>1$ are not Delzant.

The number $\nu (n)$ of all monotone Delzant polytopes of dimension~$n$ is much larger than $\nu_c (n)$: For small values of~$n$ we have
\begin{equation} \label{t:1}
\begin{array}{|c||c|c|c|c|c|c|c|c|c|c|}\hline
n      &   1&   2&   3&   4&   5&   6&   7&   8 & 9 \\ \hline
\nu_c(n) & 1&   2&   2&   4&   4&   7&   7&   12 & 12 \\ \hline
\nu(n) & 1 &  5 & 18 & 124 &  866 & 7 \, 622 & 72\,256 & 749\,892 & 8\,229\,721\\
\hline
\end{array}
\end{equation}

The next few values for $\nu_c(2k)$ are $19, 30, 45, 67, 97, 139$. The growth of $\nu_c$ is subexponential. Indeed, Since the partition function $p(n)$ grows like $e^{\sqrt{n}}$, $\nu_c(n)$ grows like $e^{\sqrt{n}}$ as well. On the other hand, for $\nu(2) =5$ see Figure~\ref{fig:2dim}. The value $\nu(3) = 18$ was found in \cite{Bat81, WatWat82}, see also \cite[pp.\ 90]{Oda88}, and $\nu (4) = 124$ was found in~\cite{Bat99, Sat00}. The values $\nu (n)$ for $5 \leq n \leq 8$ were computed by \O bro~\cite{Obr07}, and $\nu(9)$ by Paffenholz~\cite{Paf17}. The asymptotic behaviour of $\nu (n)$ is unfortunately not known, but based on discussions with Benjamin Nill and Andreas Paffenholz we expect that $\nu (n)$ grows at least exponentially. It follows that the property of a toric monotone symplectic manifold to have real central fiber is very restrictive.

\medskip
\paragraph{\bf Application II: Chekanov tori.} The Chekanov torus was defined in~\cite{Che96} as the first example of monotone Lagrangian tori in $\RR^{2n}$ which is not symplectomorphic to a product torus. We show that it can be embedded into any toric monotone symplectic manifold $M$ and compute its displacement energy germ, by closely following the ideas used in~\cite{CheSch10}. In particular, its germ shows that the Chekanov torus is exotic in~$M$. Furthermore, the polytope which is obtained as level set of the displacement energy germ is never centrally symmetric, and hence the Chekanov torus in $M$ is not real, see for example Figure~\ref{fig:level3} in Section~\ref{sec:chekanov}.

\begin{theorem}
\label{t:chek}
Let $M$ be a toric monotone symplectic manifold satisfying property $FS$. Then the Chekanov torus can be embedded into $M$ to yield an exotic Lagrangian which is not real. 
\end{theorem}
\noindent
In the case of $M= \times_n S^2$, we prove that arbitrary products of Chekanov tori are pairwise not symplectomorphic and hence we get a collection of non real exotic Lagrangian tori in $\times_n S^2$ whose cardinality grows like the partition function and hence like~$e^{\sqrt{n}}$ with $n$. In case the moment polytope of $M$ is centrally symmetric, we furthermore prove that the Chekanov torus and products thereof can be realized as the fixed point set of a smooth involution. Hence in that case, Theorem~\ref{t:chek} exhibits a symplectic phenomenon. We recall that Kim showed in \cite{Kim19b} that the Chekanov torus in $S^2 \times S^2$ is not real by using that the count of Maslov-index two holomorphic disks on this torus is five, while this count on real Lagrangians must be even. For Chekanov tori in other toric monotone symplectic manifolds, this condition for realness seems to be less useful than our condition, since the count of Maslov-index two disks is difficult, see Remark~\ref{rk:wallcrossing}.

\medskip
\paragraph{\bf Remark and Questions.}
 
In this paper we look at monotone Lagrangian tori that appear as the fibre of a torus fibration. In del~Pezzo surfaces there exist many more monotone Lagrangian tori, that are the fibre of an \emph{almost} toric fibration. An infinity of such tori were constructed by Vianna~\cite{Via17}, and for del~Pezzo surfaces different from $\CP^2$ many more (roughly an infinity for each of Vianna's tori) are found in~\cite{BreCheSch20}. Since none of the almost toric base polygons of these tori is centrally symmetric, our main result still holds, showing that none of the new tori is real. In view of this and our present work we ask:
\begin{question}
{\rm
Let $(M,\omega)$ be a toric monotone symplectic manifold which contains a real Lagrangian torus.
\begin{itemize}
\item[(1)] Is the moment polytope of $M$ necessarily centrally symmetric? 

\smallskip
\item[(2)] Suppose the moment polytope of $M$ is centrally symmetric. Is the central fibre the unique real Lagrangian, up to Hamiltonian isotopy?
\end{itemize}

\noindent
In dimension four, the only closed toric monotone symplectic manifolds 
are the toric del Pezzo surfaces $S^2 \times S^2$ and the $k$-fold blow-up~$X_k$ 
of~$\CP^2$ for $k = 0,1,2,3$,
and the only closed monotone symplectic manifolds are the del Pezzo surfaces 
$S^2 \times S^2$ and $X_k$ with $k \leq 8$, with unique symplectic structure up to scaling, 
see~\cite{Via17} for references.

\begin{itemize}
\item[(3)] Is it true that the only real tori up to Hamiltonian diffeomorphism in a 
closed monotone symplectic 4-manifold are the Clifford torus in $S^2 \times S^2$ and in $X_3$?
\end{itemize}
\noindent
We note that the uniqueness in $S^2 \times S^2$ is ongoing work by Kim \cite{Kim19c, Kim19b}.
}
\end{question} 

\medskip
\paragraph{\bf A few motivations for the study of real Lagrangian submanifolds.}

\smallskip \noindent
We conclude this introduction by mentioning some of the strands that lead to the study of real Lagrangian submanifolds.

\medskip
1. A related theme is the study of 
real algebraic varieties, namely the fixed point set of an anti-holomorphic involution of a 
complex algebraic variety.
The study of their topologocial properties has a rich history with an impressive body of results, 
see~\cite{DegKha00}.
It is interesting to see which of these results have analogues in the symplectic setting.

\smallskip
2. Let $\iota$ be a \emph{smooth} involution of a manifold $X$. Classical Smith theory  
\begin{eqnarray} 
\label{eq:smith1}
\chi (\Fix (\iota)) &=& \chi (X) \quad \mod 2, \\
\label{eq:smith2}
\dim H(\Fix (\iota);\ZZ_2) &\leq & \dim H(X;\ZZ_2)
\end{eqnarray}
relates the homology of the fixed point set of a smooth involution to the homology of the ambient manifold. We refer to \cite{Bor60} for details. It is interesting to find invariants of real Lagrangian manifolds that go beyond the Smith inequalities, and thus describe a genuine symplectic pheonomenon. As was noted by Kim~\cite{Kim19b}, the Chekanov torus in $S^2 \times S^2$ can be realized as the fixed point set of a smooth involution, but not of an anti-symplectic one. In the general context of toric monotone symplectic manifolds (see table~$(\ref{t:1})$), Theorem~\ref{t:toric} seems to yield a significantly stronger obstruction than Smith theory, which only excludes two of the five manifolds in dimension $4$ and five of the eighteen manifolds in dimension $6$.

\medskip
In symplectic geometry, real Lagrangians have appeared quite a while ago in two different forms:

\smallskip
3. Several time-honoured systems in classical mechanics, like the planar circular restricted 3-body problem, 
are invariant under several anti-symplectic involutions. Their fixed point sets can be used to find special 
orbits, see~\cite{FraVan18}.

\smallskip
4. The Arnold--Givental conjecture generalizes the classical Arnold conjecture on the number of Lagrangian intersections in terms of real Lagrangian manifolds, see~e.g.~\cite[\S 11.3 ]{McDSal17}.

\smallskip
The study of the topology of real Lagrangians in symplectic manifolds has started only very recently.
The first paper is~\cite{Kim19a}, in which Kim proved that real Lagrangians in a given compact symplectic manifold 
are unique up to cobordism, and that the only real Lagrangian in~$\CP^2$ is~$\RP^2$ up to Hamiltonian isotopy. More recent are \cite{Kim19b, Kim19c} that we discussed earlier. In collaboration with J. Kim and J. Moon~\cite{BreKimMoo19}, we construct many real Lagrangians in toric symplectic manifolds by lifting symmetries of the moment polytope. 

\medskip
\paragraph{\bf Organization of the paper.}
In Section~\ref{sec:versal} we discuss real Lagrangians and versal deformations. We prove Theorem \ref{t:dis} on the displacement energy germ of real Lagrangians. In Section~\ref{sec:displacement}, we discuss the displacement energy of toric fibres with a focus on the case in which the moment polytope has property $FS$. This discussion is instrumental for both our applications. In Section~\ref{sec:toric}, we discuss whether fibres of toric symplectic manifolds are real and establish a criterion in terms of the geometry of the corresponding moment polytope. In particular, we prove Theorem~\ref{t:toric}. In Section~\ref{sec:chekanov}, we deal with Chekanov tori in toric monotone symplectic manifolds and show that none of them are real, see Theorem~\ref{t:chek}. The appendix in Section~\ref{sec:appendix} outlines an alternate approach to our results using~$J$-holomorphic curves.

\medskip
\paragraph{\bf Acknowledgements.} I thank Joontae Kim for introducing me to real symplectic geometry and for his invitation to KIAS in May 2019, where the main idea of this paper arose from numerous stimulating discussions. I also wish to thank Yuri Chekanov for useful remarks and for agreeing that I use some of our joint results from~\cite{BreCheSch20} in this paper. Many thanks to Andreas Paffenholz for having written a script checking
that all of the Delzant polytopes in dimension~$9$ have the required property $FS$ and to Grigory Mikhalkin and an anonymous referee for suggesting the approach to the problem outlined in the appendix. I am grateful to Felix Schlenk for his generous support and countless remarks from which this paper has greatly benefited.

\section{Versal Deformations of real Lagrangians}
\label{sec:versal}
In this section, we will discuss real Lagrangians, displacement energy and versal deformations. In particular, we will prove Theorem~\ref{t:dis}, the proof of which relies on two key observations. Firstly, the displacement energy is invariant under anti-symplectic involutions, see Proposition~\ref{prop:einvariance}. Secondly, if we combine this invariance with a $\ZZ_2$-equivariant Weinstein neighbourhood Theorem, we obtain the desired result.

\subsection{Real Lagrangians} 
Let $(M,\omega)$ be a symplectic manifold and let~$\sigma$ be an anti-symplectic involution on $M$, i.e.\ a smooth map satisfying~$\sigma \circ \sigma = \idd$ and $\sigma^*\omega = -\omega$. Its fixed point set $\Fix \sigma$ is a (possibly not connected) Lagrangian submanifold whenever it is not empty. 

\begin{definition}
A Lagrangian submanifold $L$ in $(M,\omega)$ is called \emph{real} if there is an anti-symplectic involution of $M$ having $L$ as a connected component of its fixed point set. 
\end{definition}

\begin{example}
\label{ex:equator}
{\rm
The equator in the standard symplectic $2$-sphere $(S^2,\omega)$ is real. The corresponding involution is given by reversing the height $z\mapsto -z$. By taking the product of this example, we can describe the product of equators (also known as the Clifford torus) as the fixed point set of an anti-symplectic involution on $\times_n S^2$. 
}
\end{example}

\begin{example}
{\rm
Let $(\CP^n, \omega_{\rm FS})$ be the complex projective space equipped with the Fubini--Study form. Then $\RP^n \subset \CP^n$ is real since it is the fixed point set of the anti-symplectic involution
	\begin{equation*}
		\sigma \colon \CP^n \rightarrow \CP^n, \quad [z_0 : \ldots : z_n] \mapsto [\overline{z}_0: \ldots :\overline{z}_n].
	\end{equation*}
\noindent 
It is well-known that this example can be generalized to all toric manifolds, an observation which gives rise to so-called real toric geometry. See for example~\cite{DegKha00}.
}
\end{example}

\begin{example}
\label{ex:antidiag}
{\rm Any symplectic manifold $(M,\omega)$ can be seen as a real Lagrangian submanifold in $(M\times M, \omega \oplus - \omega)$. The embedding is given by the diagonal map $p \mapsto (p,p)$ and the corresponding anti-symplectic involution is given by exchanging the two coordinates in $M\times M$.
}
\end{example}

\begin{example} 
\label{ex:cotangent}
{\rm
Let $(T^*Q,\omega_0 = -d\lambda)$ be the cotangent bundle of a smooth manifold $Q$ equipped with its canonical symplectic form. The map which reverses momenta,
\begin{equation*} 
	\sigma_0 \colon T^*Q \rightarrow T^*Q, \quad (q,p) \mapsto (q,-p),
\end{equation*}
satisfies $\sigma_0^*\lambda = -\lambda$ and is therefore an anti-symplectic involution. Its fixed point set is the zero section
\begin{equation*} 
	\Fix \sigma_0 = Q \subset T^*Q.
\end{equation*}
}
\end{example}

By Weinstein's Lagrangian neighbourhood theorem and Example~\ref{ex:cotangent}, any Lagrangian submanifold admits a \textit{locally defined} anti-symplectic involution of which it is the fixed point set. Of course, locally defined involutions might not extend globally. On the other hand, Meyer~\cite{Mey81} proved that any anti-symplectic involution $\sigma$ with non-empty fixed point set is locally of the form described in Example~\ref{ex:cotangent}. This can be viewed as a $\mathbb{Z}_2$-equivariant version of Weinstein's theorem.

\begin{theorem} 
{\bf (Meyer~\cite{Mey81})}
\label{thm:equiweinstein}
Let $\sigma$ be an anti-symplectic involution of a symplectic manifold $(M,\omega)$ containing a Lagrangian $L\subseteq \Fix \sigma \neq \varnothing$. Furthermore let $T^*L$ be equipped with its canonical symplectic form and the anti-symplectic involution $\sigma_0$ which reverses momenta. Then there is a $\sigma$-invariant neighbourhood $V$ of $L$, a $\sigma_0$-invariant neighbourhood~$U$ of the zero-section in $T^*L$ and a symplectomorphism
	\begin{equation*}
		g \colon (U, \omega_0\vert_U) \rightarrow (V, \omega\vert_V),
	\end{equation*}
which maps the zero section to $L$ and which intertwines the anti-symplectic involutions $\sigma$ and $\sigma_0$,
	\begin{equation}
	\label{eq:equiweinstein}
		g \circ \sigma_0 = \sigma \circ g.
	\end{equation}
\end{theorem}

\subsection{Displacement energy}
Recall that the displacement energy of a compact subset $A$ of a symplectic manifold $(M,\omega)$ is defined as
\begin{equation*}
	e_M(A)= \inf \left\lbrace \Vert H \Vert \left\vert \; H \in C_{c}^{\infty}([0,1] \times M),\, \varphi_H^1 (A) \cap A = \varnothing \right.\right\rbrace,
\end{equation*}

\noindent
where
\begin{equation*}
\Vert H \Vert = \int_0^1 \left( \max_{x\in M} H_t(x) - \min_{x \in M} H_t(x) \right) dt
\end{equation*}
is the Hofer norm on $C_{c}^{\infty}([0,1] \times M)$. By convention, we put $e_M(A)=\infty$ whenever the set of displacements is empty.

\begin{example}
\label{ex:productenergy}
{ \rm
Let $T(a) \subset (\RR^2,\omega_0)$ be the circle enclosing area $a > 0$ in the plane. Its displacement energy is 
\begin{equation*}
	e_{\RR^2}(T(a))=a.
\end{equation*}
By taking products, we obtain Lagrangian \emph{product tori} $T(a_1,\dots,a_n) = T(a_1) \times \dots \times T(a_n) \subset (\RR^{2n}, \omega_0)$. Their displacement energy is (see Remark~\ref{rk:productenergy})
\begin{equation*}
	e_{\RR^{2n}}(T(a_1,\dots,a_n)) = \min\{a_1,\dots,a_n\}.
\end{equation*}
}
\end{example}

Given a symplectomorphism $\psi$ of $(M,\omega)$ we have $\varphi_{H \circ \psi^{-1}}^t = \psi \circ \varphi_H^t \circ \psi^{-1}$. The set $A$ is thus displaced by the time-one map of $H$ if and only if $\psi(A)$ is displaced by the time-one map of $H \circ \psi^{-1}$. Since the Hofer norm satisfies $\Vert H \circ \psi^{-1} \Vert = \Vert H \Vert$, it follows that displacement energy is invariant under symplectomorphisms,
\begin{equation*} 
	e_M(\psi (A)) = e_M(A).
\end{equation*}
The same is true for anti-symplectic involutions. 

\begin{proposition}
\label{prop:einvariance}
Let $\sigma$ be an anti-symplectic involution on a symplectic manifold $(M,\omega)$. Then the displacement energy is invariant under $\sigma$ in the sense that 
\begin{equation*}
	e_M(\sigma(A))=e_M(A)
\end{equation*}
for any compact subset $A \subset M$.
\end{proposition}

\proof Let $H \in C_{c}^{\infty}([0,1] \times M)$ be a Hamiltonian, $X_H^t$ and $\varphi_H^t$ its associated vector field and flow. Since $\sigma$ is an anti-symplectic involution, we have 
\begin{equation*}
	X^t_{H \circ \sigma} = -\sigma_* (X^t_H \circ \sigma).
\end{equation*}
Define $H'_t = -H_t \circ \sigma$. Its Hamiltonian vector field is
\begin{equation*}
	X^t_{H'} = \sigma_* (X^t_H \circ \sigma) 
\end{equation*}
and thus we get for the respective flows
\begin{equation*}
\varphi_{H'}^t = \sigma \circ \varphi_H^{t} \circ \sigma.
\end{equation*}
This proves that a set $A$ is displaced by $\varphi_H^1$ if and only if $\sigma(A)$ is displaced by $\varphi_{H'}^1$. Since $\Vert H' \Vert = \Vert H \Vert$, the claim follows.
\proofend

\subsection{Versal Deformations}

Versal deformations were introduced in~\cite{Che96} and subsequently used in~\cite{CheSch10} and~\cite{CheSch16} as a tool to distinguish Lagrangian submanifolds. The idea is to look at the behaviour of known symplectic invariants on neighbouring Lagrangians of the submanifolds in question. Let us outline the construction. Since we will only use the displacement energy as an invariant, we will restrict ourselves to this case. We refer to~\cite{CheSch16} for details.\\

In every cotangent bundle $T^*L$ of a closed Lagrangian submanifold, Lagrangians which are $C^1$-close to the zero section can be identified with the graphs of closed one-forms. Using Weinstein's theorem, one can translate this identification to the case of any Lagrangian $L \subset (M,\omega)$ as follows. For a given Weinstein chart $g : T^*L \supset U \rightarrow V \subset M $ there is a $C^1$-neighbourhood $\widehat{\mathcal{U}} \subset \Omega^1_{{\rm cl}}(L)$ of the zero section in the space of closed one-forms, a $C^1$-neighbourhood $\widehat{\mathcal{V}}$ of $L$ in the space of Lagrangian submanifolds in $M$, and a bijection
\begin{equation*}
	\widehat{w}^g_L \colon \widehat{\mathcal{U}} \rightarrow \widehat{\mathcal{V}}, \quad \alpha \mapsto g(\Gamma_{\alpha}),
\end{equation*}
where we denote the graph of $\alpha \in \Omega^1(L)$ by $\Gamma_{\alpha}$. Furthermore, $C^1$-small Hamiltonian perturbations of the zero section in $T^*L$ are in one-to-one correspondence with $C^1$-small exact one-forms, and hence the above map descends to 
\begin{equation*}
	w^g_L \colon \mathcal{U} \rightarrow \mathcal{V},
\end{equation*}
where we divide out exact one-forms on the left-hand side and Hamiltonian isotopies on the right-hand side. In particular we can view $\mathcal{U}$ as a neighbourhood of zero in $H^1(L,\mathbb{R})$. Up to Hamiltonian isotopy, neighbouring Lagrangians of $L$ are thus parametrized by a neighbourhood of zero in the vector space $H^1(L,\mathbb{R})$.\\

As displacement energy is invariant under Hamiltonian isotopies we can compose it with the above map $w_L^g$ to obtain a function on $\mathcal{U}$
\begin{equation*}
	 H^1(L,\mathbb{R}) \supset \mathcal{U} \rightarrow \mathbb{R} \cup \{\infty\}, \quad [\alpha ] \mapsto e_M(g(\Gamma_{\alpha})).
\end{equation*}
The germ at $0$ associated to this function corresponds to the displacement energy of neighbouring Lagrangians of $L$ and will be denoted by 
\begin{equation*}
	 S_L^g \colon (H^1(L,\mathbb{R}),0) \rightarrow \mathbb{R} \cup \{\infty\}.
\end{equation*}
\noindent
The following remark is crucial for what will follow. 

\begin{remark}
\label{rk:weinstein_indep}
{\rm 
The germ of the bijection $w^g_L$ is independent of the choice of Weinstein chart $g$ and thus so is the germ $S^g_L$. Hence we will write $S_L = S^g_L$. See~\cite{CheSch16} for details.
}
\end{remark}
\noindent
We are now in a position to prove Theorem \ref{t:dis}, which we recall for the reader's convenience.

\begin{theorem} 
\label{thm:main}
Assume that $L \subseteq \Fix \sigma$ is a compact real Lagrangian submanifold of $(M,\omega)$. 
Then the displacement energy germ $S_L$ is even,
\begin{equation*}
	S_L(-p) = S_L(p).
\end{equation*}
\end{theorem}
\proof 
By Theorem~\ref{thm:equiweinstein} we can pick a Weinstein neighbourhood $g$ such that $g \circ \sigma_0 = \sigma \circ g$. Let $\alpha \in \Omega^1_{\rm cl}(L)$ be a one-form representing $p$, then $\sigma_0(\Gamma_{\alpha}) = \Gamma_{-\alpha}$. Hence, using the invariance of the displacement energy under anti-symplectic involutions, we find
	\begin{eqnarray*}
		S^g_L (-\alpha)
		&=& e_M(g(\Gamma_{-\alpha})) \\
		&=& e_M(g ( \sigma_0 (\Gamma_{\alpha}))) \\
		&=& e_M(\sigma(g(\Gamma_{\alpha}))) \\
		&=& e_M(g(\Gamma_{\alpha})) \\
		&=& S^g_L(\alpha).
	\end{eqnarray*}
Since $S^g_L = S_L$ is independent of the choice of $g$, the claim follows. \proofend

\section{Displacement energy of toric fibres}
\label{sec:displacement}

In this section we compute the displacement energy of toric fibres. We begin by proving that displacement energy can only increase under symplectic reduction. This observation was already made in \cite{AbrMac13} and will be used here to prove the existence of a lower bound as well as an upper bound on the displacement energy of toric fibres. For the lower bound, we will use the fact that any toric symplectic manifold can be seen as a symplectic quotient of some $\CC^k$ via Delzant's construction. For the upper bound, we will give a slightly modified version of McDuff's method by probes, see \cite{McD11}. In the last part of this section we will apply these results to compute the displacement energy of toric fibres in toric monotone symplectic manifolds. This is a crucial ingredient for Sections~\ref{sec:toric} and~\ref{sec:chekanov}.

\subsection{Displacement energy and symplectic reduction.} 
Let $(\widehat{M} , \widehat{\omega} , \nu)$ be a Hamiltonian $G$-space which admits symplectic reduction at the level $0 \in \mathfrak{g}^*$, i.e. $0$ is a regular value and $G$ acts freely on $Z = \nu^{-1}(0)$. This means that we have the following reduction diagram

\begin{center}
	\begin{tikzcd}
		Z = \nu^{-1}(0) \arrow[hook]{r} \arrow[two heads]{d}{p}
		& (\widehat{M},\widehat{\omega}) \\
		\momega
	\end{tikzcd}
\end{center}

\noindent
with $\widehat{\omega}\vert_{TZ} = p^*\omega$. Furthermore, assume that the symplectic quotient $\momega$ is compact. 

\begin{lemma}
\label{lem:energyreduction}
	Under the above hypotheses, we have 
	\begin{equation*}
		e_{\widehat{M}}(p^{-1}(A))	\leqslant e_M(A)
	\end{equation*}
	for any set $A \subset M$.
\end{lemma}

\noindent
In other words, symplectic reduction can only increase displacement energy. The proof of Lemma \ref{lem:energyreduction} runs as follows. For any Hamiltonian $H \in C^{\infty}(M \times [0,1])$ which displaces $A$ we will construct a compactly supported Hamiltonian $\widehat{H} \in C_c^{\infty}(\widehat{M} \times [0,1])$ which displaces $p^{-1}(A)$ and which has the same Hofer norm as $H$. The Hamiltonian $\widehat{H}$ is obtained as an extension of the lift $p^*H \in C^{\infty}(Z \times [0,1])$ which is zero outside of a tubular neighbourhood of $Z \subset \widehat{M}$. Although this was already outlined in \cite{AbrMac13}, we give a full proof for the reader's convenience. \\

\proofof{Lemma \ref{lem:energyreduction}} If $A$ is not displaceable, there is nothing to show. Therefore let $H$ be a Hamiltonian on $M$ which displaces $A$. We can assume that $\min_{p\in M}H_t(p) = 0$ for all $t\in [0,1]$. Now fix a time~$t\in [0,1]$ and pick a tubular neighbourhood of $Z$, i.e. a diffeomorphism 
	\begin{equation*}
		\chi : NZ \supset U \rightarrow V \subset \widehat{M}
	\end{equation*}
from a neighbourhood $U$ of the zero section inside the normal bundle $\pi : NZ \rightarrow Z$ to a neighbourhood $V$ of $Z \subset \widehat{M}$ mapping the zero section to $Z$. Let $\rho \in C^{\infty}(U)$ be a function such that 
\begin{enumerate}
	\item[1.] $\rho = 1$ on the zero section and $\rho \leqslant 1$ elsewhere,
	\item[2.] $\rho$ is compactly supported.
\end{enumerate} 
We can now define $\widehat{H}_t$ on $U$ by putting $\widehat{H}_t(v) = \rho(v)p^*H_t(\pi(v))$. By using $\chi$, we transport this function to a function $\widehat{H}_t$ on $V$, which can be smoothly extended to all of $\widehat{M}$ by zero since $\rho$ has compact support. Notice that the Hofer norm of $\widehat{H}$ is equal to the Hofer norm of $H$. For $\widehat{H}_t \in C^{\infty}(\widehat{M})$ we have
	\begin{equation}
		\label{eq:hbar}
		\widehat{H}_t\vert_Z = p^*H_t.
	\end{equation}
In particular, $\widehat{H}_t\vert_Z$ is invariant under the $G$-action on $Z$. We will show that the restriction of the Hamiltonian vector field $X_{\widehat{H}}^t$ to $Z$ 
	\begin{enumerate}
		\item[1.] is tangent to $Z$,
			\begin{equation}
				\label{eq:tangency}
				(X_{\widehat{H}}^t)_z \in T_zZ \quad \forall z \in Z;
			\end{equation}
		\item[2.] projects to the Hamiltonian vector field of $H$ on $M$,
			\begin{equation}
				\label{eq:projection}
				p_* (X_{\widehat{H}}^t\vert_Z) = X_H^t.
			\end{equation}
	\end{enumerate}
\noindent
In order to prove $(\ref{eq:tangency})$, we use the invariance of $\widehat{H}_t\vert_Z$ under the action of $G$, which implies that the following equivalent conditions hold
	\begin{eqnarray*}
		&& d\widehat{H}_t(z) (X_{\zeta})_z = 0  \quad \forall \zeta \in \mathfrak{g}, \\	
		&\Leftrightarrow & \langle d\nu(z)(X^t_{\widehat{H}})_z , \zeta \rangle = 0 \quad \forall \zeta \in \mathfrak{g},\\
		&\Leftrightarrow & (X^t_{\widehat{H}})_z \in T_zZ.
	\end{eqnarray*}
The last line follows from the fact that $T_zZ = T_z\nu^{-1}(0)=\ker d\nu(z)$. Let $Y \in TM$ and pick $\widehat{Y} \in TZ$ so that $p_*\widehat{Y} = Y$. Using $(\ref{eq:hbar})$, we find $d(\widehat{H}_t\vert_Z)(\widehat{Y}) = dH_t(Y)$, which we use to compute
	\begin{eqnarray*}
		\omega(p_*X_{\widehat{H}}^t , Y)
		&=& \omega(p_*X_{\widehat{H}}^t , p_*\widehat{Y}) \\
		&=& (p^*\omega)(X_{\widehat{H}}^t , \widehat{Y}) \\
		&=& \widehat{\omega}(X_{\widehat{H}}^t , \widehat{Y}) \\
		&=& d(\widehat{H}_t\vert_Z)(\widehat{Y}) \\
		&=& dH_t(Y) \\
		&=& \omega(X_H^t,Y).
	\end{eqnarray*}
This proves $(\ref{eq:projection})$. Now let $\varphi_H^t$ and $\varphi_{\widehat{H}}^t$ denote the corresponding Hamiltonian flows. Since equation $(\ref{eq:projection})$ holds for all $t\in [0,1]$, we have
	\begin{equation}
		\label{eq:projection_flow}
		p \circ \varphi_{\widehat{H}}^t\vert_Z = \varphi_H^t \circ p.
	\end{equation}
Since $\varphi_H^1(A) \cap A = \varnothing$, take the pre-image under $p$ of both sides to get $p^{-1}(\varphi_H^1(A)) \cap p^{-1}(A) = \varnothing$. Together with equation $(\ref{eq:projection_flow})$,
\begin{equation*}
	\varphi_{\widehat{H}}^1(p^{-1}(A)) \cap p^{-1}(A) \subseteq p^{-1}(\varphi_H^1(A)) \cap p^{-1}(A) = \varnothing 
\end{equation*}
and hence $\varphi_{\widehat{H}}^1$ displaces $p^{-1}(A)$.
\proofend

\subsection{Lower bound for toric fibres} \label{ss:lbtoric} Let $(M^{2n},\omega)$ be a compact toric symplectic manifold. By this we mean that $T^n$ acts effectively on $M$ by Hamiltonian diffeomorphisms which are generated by a moment map $\mu \colon M \rightarrow \mathfrak{t}^*$. We identify the dual of the Lie algebra of $T^n$ with $\RR^n$ by choice of a basis. As is the case for all Hamiltonian torus actions, the image of $\mu$ is a convex polytope $\Delta = \mu(M) \subset \RR^n$, called \emph{moment polytope}. Since $M$ is toric, the corresponding moment polytope has the Delzant property, see \cite{Del88} or \cite{Aud04} for details. Furthermore, Delzant showed that $M$ can be reconstructed from such $\Delta$ by taking a suitable symplectic quotient of $\CC^k$ by the action of a linear subtorus of $T^k$ acting by the standard action on $\CC^k$. Let $\nu$ be the moment map of this action. The situation is summarized by the following reduction diagram

\begin{center}
	\begin{tikzcd}
		Z = \nu^{-1}(0) \arrow[hook]{r} \arrow[two heads]{d}{p}
		& (\CC^k,\omega_0) \\
		\momega \arrow{r}{\mu} 
		& \Delta.
	\end{tikzcd}
\end{center}

\noindent 
We describe the moment polytope $\Delta \subset \RR^n$ of $M$ by a set of inequalities 
\begin{equation*}
	\langle x , v_i \rangle \leqslant \kappa_i , \quad i \in \{1,\dots,k\},
\end{equation*}
where the $v_i$ are the unique outward-pointing normal vectors to the facets of $\Delta$ which are primitive in the lattice $\ZZ^n \subset \RR^n$. Define the functionals on $\RR^n$
\begin{equation*}
\ell_i(x) = \kappa_i - \langle x , v_i \rangle 
\end{equation*}
for all $i \in \{1,\dots,k\}$. Every $\ell_i$ defines a half-space $\{\ell_i \geq 0\}$ and the moment polytope $\Delta$ is given by the intersection of these half-spaces. Using Lemma \ref{lem:energyreduction}, we will give a lower bound for the displacement energy of any toric fibre $T_x = \mu^{-1}(x)$.

\begin{proposition}
\label{prop:lower}
Let $\momegam$ be a toric symplectic manifold with moment polytope $\Delta$. Then for every $x \in \Delta$ the displacement energy of the corresponding toric fibre is bounded from below by 
	\begin{equation*}
		e_{\Delta}(x) = e_M(T_x) \geqslant \min \{\ell_1(x),\ldots,\ell_k(x)\},
	\end{equation*}
where $\ell_i(x)$ is the affine distance of $x$ to the $i$-th facet of $\Delta$.
\end{proposition}

\proof
As is clear from the Delzant construction, $p^{-1}(T_x)=p^{-1}(\mu^{-1}(x)) \subset \CC^k$ is the product torus
\begin{equation*}
	T(a_1,...,a_k) = \{(z_1,...,z_k) \in \CC^k \;\vert\; \pi\vert z_i \vert^2 = a_i \},
\end{equation*}
with $a_i = \ell_i(x)$. Since $e_{\CC^k}(T(a_1,...,a_k)) = \min \{a_1,\dots,a_k\}$ by Remark~\ref{rk:productenergy}, the claim follows from Lemma \ref{lem:energyreduction}.
\proofend

\begin{remark}
\label{rk:productenergy}
{\rm
In order to compute the displacement energy of a product torus in~$\CC^n$, we use the inequalities
\begin{equation*}
\min \{a_1,\dots,a_k\} \leq c_1(T(a_1,...,a_k)) \leq e_{\CC^k}(T(a_1,...,a_k)) \leq \min \{a_1,\dots,a_k\},
\end{equation*}
where~$c_1$ denotes the first Ekeland-Hofer capacity. The first inequality follows from Theorem~$(b)$ on page~43 of \cite{Sik90} and the second from~\cite[Theorem 1.6]{Hof90}, which are both obtained by applying the calculus of variations to the action functional of classical mechanics. The third inequality follows from Proposition~\ref{prop:upper}. This is the only \emph{hard} symplectic result we use and hence our methods do not rely on~$J$-holomorphic curves, with the obvious exception of the complementary appendix. 
}
\end{remark}

\subsection{Upper bound for toric fibres}
\begin{figure}[h]
		\begin{tikzpicture}[scale=1.3]
			\draw [thin,fill=black!15] (-1,0)--(3,-1)--(2,2)--(-1,0);
			\draw[thick] (-1,0)--(3,-1)--(-1,0);
			\draw[thick,dashed]  (1,-0.5)--(2.5,0.5)--(1,-0.5);
			\draw [thick,black,->] (1,-0.5)--(0,-1.16);
			\fill[thick, black] (1,-0.5)  circle[radius=1.5pt];
			\fill[thick, black] (2.5,0.5)  circle[radius=1.5pt];
			\fill[thick, black!0] (2.5,0.5)  circle[radius=1pt];
			\fill[thick, black] (1.5,-0.16)  circle[radius=1.5pt];
			\node at (0.2,-1.25){$u$};
			\node at (1,-0.75){$w$};
			\node at (1.5,-0.4){$x$};
			\node at (1.7,0.45){$P_{i,u}(w)$};
			\node at (-0.25,-0.45){$F_i$};
		\end{tikzpicture}
	\caption{The probe $P_{i,u}(w)$.}
	\label{fig:probe}
\end{figure}
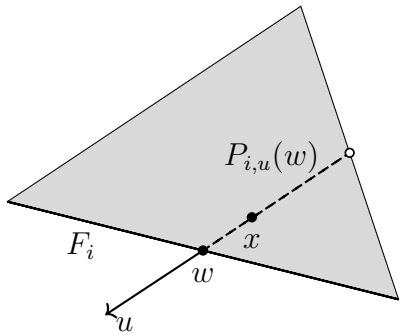

In order to prove displaceability in toric symplectic manifolds, \mbox{McDuff} introduced probes in \cite{McD11}, a technique independently found in \cite{CheSch10}. We will show that probes can be interpreted in the framework of Lemma \ref{lem:energyreduction} by performing symplectic reduction on the pre-image of the probe. Let $\momega$ be a toric symplectic manifold with moment map $\mu$ and moment polytope $\Delta = \{\ell_i \geqslant 0, \; \forall i\}$. A probe $P_{i,u}(w)$ is determined by a facet $F_i=\{\ell_i = 0\} \cap \Delta$ of $\Delta$, a point $w \in F_i$ and a vector $u \in \ZZ^n$ which is integrally transverse to $F_i$. By this we mean that $u$ can be completed to a $\ZZ$-basis of $\ZZ^n$ by vectors parallel to $F_i$. The set $P_{i,u}(w) \subset \RR^n$ is the half open line segment obtained as the union of $\{w\}$ with the open segment defined by the intersection of $\mathring{\Delta}$ with the line emanating from $w$ in direction $u$, see Figure \ref{fig:probe}. Displaceability of toric fibres lying on a suitable probe was proved in \cite{McD11}. 

\begin{proposition}
\label{prop:upper}
Let $x\in \Delta$ be a point in a probe $P_{i,u}(w)$ lying in the same half of $P_{i,u}(w)$ as $w$ and not on the midpoint of the probe. Then 
	\begin{equation*}
		e_{\Delta}(x) \leqslant \ell_i(x).
	\end{equation*}
\end{proposition}

\proof 
Since $u$ is integrally transverse to $F_i$ we can assume, up to applying a transformation in $\GL(n,\ZZ)$, that $u = e_1$ and that $F_i$ lies in the hyperplane spanned by $e_2,\dots,e_n$. Hence $w= (0,w')$ for some $w' \in \RR^{n-1}$ and $x=(\ell_i(x),w')$. Let $U = \mu^{-1}(\mathring{\Delta} \cup \mathring{F_i}) \subset M$. The subtorus $T^{n-1} = \{1\} \times S^1 \times \dots \times S^1 \subset T^n$ acts freely\footnote{This can be seen as follows. For toric manifolds $\momegam$ the stabilizer of any point $p\in M$ can be read off from the moment polytope (viewed as $\Delta \subset \mathfrak{t}^*$) by taking the subtorus inside $T^n$ which is generated by the annihilator of the smallest face of $\Delta$ which contains $\mu(p)$. In our situation, the annihilator is generated by $e_1$ and hence the stabilizer is the first coordinate circle in $T^n$.} on $U$. The moment map of this action $\mu' : U \rightarrow \RR^{n-1}$ is obtained by restricting $\mu$ to $U$ and by dropping the first coordinate 
\begin{equation*}
	\mu'(y_1,\dots,y_{n-1}) = (\mu_2\vert_{U}(y_1),\dots,\mu_n\vert_{U}(y_{n-1})).
\end{equation*}
We get $(\mu')^{-1}(w') = \mu^{-1}(P_{i,u}(w))$ and since $T^{n-1}$ acts freely on this set, we can consider the following symplectic reduction

\begin{center}
	\begin{tikzcd}
		(\mu')^{-1}(w') \arrow[hook]{r} \arrow[two heads]{d}{p}
		& (U,\omega\vert_U) \\
		(D^2(a),\omega_0).
	\end{tikzcd}
\end{center}
Here, the reduced space is an open disk of area $a$ equal to the affine length of the probe. The fibre we are interested in is 
\begin{equation*}
	T_x= \mu^{-1}(x) = p^{-1}(T(\ell_i(x))),
\end{equation*}
where $T(\ell_i(x)) \subset D^2(a)$ is the circle bounding area $\ell_i(x)$. By our assumption on $x$, we have $\ell_i(x) < \frac{a}{2}$ and therefore $T(\ell_i(x)) \subset D^2(a)$ has displacement energy $\ell_i(x)$. Hence by Lemma \ref{lem:energyreduction} and the fact that $U \subset M$, we find 
\begin{equation*}
	e_{\Delta}(x) = e_M(T_x) \leqslant e_{U}(T_x) \leqslant e_{D^2(a)} (T(\ell_i(x))) = \ell_i(x).
\end{equation*}

\proofend

\subsection{Probes in monotone polytopes.} \label{ssec:monotone} Let $\Delta = \{\ell_i \geqslant 0, \forall i\}$ be the moment polytope of a toric monotone symplectic manifold $(M,\omega)$. We can assume that the barycentre lies in $0 \in \RR^n$ and that $\ell_i(0)=\kappa_i=1$. McDuff \cite{McD11} discovered that displaceability by probes is related to the Ewald conjecture and proved that every point except the barycentre is displaceable by probes if and only if $\Delta$ satisfies the \textit{star Ewald condition}. Since we only need to know the function $e_{\Delta} \colon x \mapsto e_M(T_x)$ on an open and dense subset of $\Delta$, we can work directly with a variation of the Ewald conjecture which has been checked by \O bro \cite{Obr07} for dimensions $\leqslant 8$ and by Paffenholz \cite{Paf17} for dimension $9$. This approach is also used in \cite{BreCheSch20}.\\

Let $\mathcal{S}(\Delta)= \Delta \cap (-\Delta) \cap \ZZ^n \setminus \{0\}$ be the set of non-zero symmetric integral points of $\Delta$.

\begin{definition}
The polytope $\Delta$ has property $FS$ if every facet $F \subset \Delta$ contains a point of the set $\mathcal{S}(\Delta)$.
\end{definition}

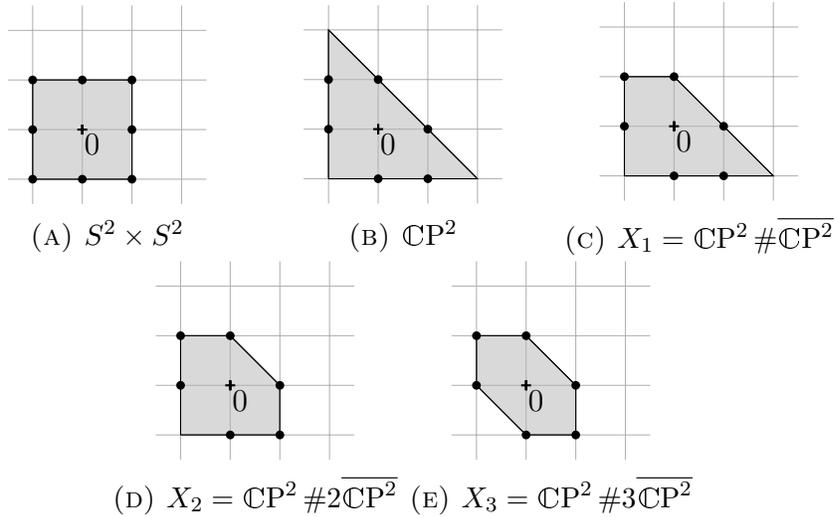
\begin{figure}[h]
	\begin{subfigure}{0.3\textwidth}
  		\centering
		\begin{tikzpicture}[scale=0.33]
			\draw [thin,fill=black!15] (-2,-2)--(2,-2)--(2,2)--(-2,2)--(-2,-2);
			\draw[step=2.0,black!30,thin] (-3,-3) grid (5,5);
			\draw[thin] (-2,-2)--(2,-2)--(2,2)--(-2,2)--(-2,-2);
			\draw[thick] (-0.2,0)--(0.2,0)--(0.2,0);
			\draw[thick] (0,0.2)--(0,-0.2)--(0,0.2);
			\fill[thick, black] (2,0)  circle[radius=5pt];
			\fill[thick, black] (2,-2)  circle[radius=5pt];
			\fill[thick, black] (0,2)  circle[radius=5pt];
			\fill[thick, black] (0,-2)  circle[radius=5pt];
			\fill[thick, black] (-2,0)  circle[radius=5pt];
			\fill[thick, black] (-2,2)  circle[radius=5pt];
			\fill[thick, black] (2,2)  circle[radius=5pt];
			\fill[thick, black] (-2,-2)  circle[radius=5pt];
			\node at (0.4,-0.6){$0$};
		\end{tikzpicture}
  	\caption{$S^2\times S^2$}
	\end{subfigure}
	\begin{subfigure}{0.3\textwidth}
		\centering
		\begin{tikzpicture}[scale=0.33]
			\draw [thin,fill=black!15] (-2,-2)--(4,-2)--(-2,4)--(-2,-2);
			\draw[step=2.0,black!30,thin] (-3,-3) grid (5,5);
			\draw[thin]  (-2,-2)--(4,-2)--(-2,4)--(-2,-2);
			\draw[thick] (-0.2,0)--(0.2,0)--(0.2,0);
			\draw[thick] (0,0.2)--(0,-0.2)--(0,0.2);
			\fill[thick, black] (2,0)  circle[radius=5pt];
			\fill[thick, black] (2,-2)  circle[radius=5pt];
			\fill[thick, black] (0,2)  circle[radius=5pt];
			\fill[thick, black] (0,-2)  circle[radius=5pt];
			\fill[thick, black] (-2,0)  circle[radius=5pt];
			\fill[thick, black] (-2,2)  circle[radius=5pt];
			\node at (0.4,-0.6){$0$};
		\end{tikzpicture}
	\caption{$\CP^2$}
	\end{subfigure}
	\begin{subfigure}{0.3\textwidth}
  		\centering
		\begin{tikzpicture}[scale=0.33]
			\draw [thin,fill=black!15] (-2,-2)--(4,-2)--(0,2)--(-2,2)--(-2,-2);
			\draw[step=2.0,black!30,thin] (-3,-3) grid (5,5);
			\draw[thin]  (-2,-2)--(4,-2)--(0,2)--(-2,2)--(-2,-2);
			\draw[thick] (-0.2,0)--(0.2,0)--(0.2,0);
			\draw[thick] (0,0.2)--(0,-0.2)--(0,0.2);
			\fill[thick, black] (2,0)  circle[radius=5pt];
			\fill[thick, black] (2,-2)  circle[radius=5pt];
			\fill[thick, black] (0,2)  circle[radius=5pt];
			\fill[thick, black] (0,-2)  circle[radius=5pt];
			\fill[thick, black] (-2,0)  circle[radius=5pt];
			\fill[thick, black] (-2,2)  circle[radius=5pt];
			\node at (0.4,-0.6){$0$};
		\end{tikzpicture}
  	\caption{$X_1 = \CP^2 \# \overline{\CP^2}$}
	\end{subfigure}
	
	\begin{subfigure}{0.3\textwidth}
  		\centering
		\begin{tikzpicture}[scale=0.33]
			\draw [thin,fill=black!15] (-2,-2)--(2,-2)--(2,0)--(0,2)--(-2,2)--(-2,-2);
			\draw[step=2.0,black!30,thin] (-3,-3) grid (5,5);
			\draw[thin]  (-2,-2)--(2,-2)--(2,0)--(0,2)--(-2,2)--(-2,-2);
			\draw[thick] (-0.2,0)--(0.2,0)--(0.2,0);
			\draw[thick] (0,0.2)--(0,-0.2)--(0,0.2);
			\fill[thick, black] (2,0)  circle[radius=5pt];
			\fill[thick, black] (2,-2)  circle[radius=5pt];
			\fill[thick, black] (0,2)  circle[radius=5pt];
			\fill[thick, black] (0,-2)  circle[radius=5pt];
			\fill[thick, black] (-2,0)  circle[radius=5pt];
			\fill[thick, black] (-2,2)  circle[radius=5pt];
			\node at (0.4,-0.6){$0$};
		\end{tikzpicture}
  	\caption{$X_2 = \CP^2 \# 2\overline{\CP^2}$}
	\end{subfigure}
	\begin{subfigure}{0.3\textwidth}
  		\centering
		\begin{tikzpicture}[scale=0.33]
			\draw [thin,fill=black!15] (-2,0)--(0,-2)--(2,-2)--(2,0)--(0,2)--(-2,2)--(-2,0);
			\draw[step=2.0,black!30,thin] (-3,-3) grid (5,5);
			\draw[thin] (-2,0)--(0,-2)--(2,-2)--(2,0)--(0,2)--(-2,2)--(-2,0);
			\draw[thick] (-0.2,0)--(0.2,0)--(0.2,0);
			\draw[thick] (0,0.2)--(0,-0.2)--(0,0.2);
			\fill[thick, black] (2,0)  circle[radius=5pt];
			\fill[thick, black] (2,-2)  circle[radius=5pt];
			\fill[thick, black] (0,2)  circle[radius=5pt];
			\fill[thick, black] (0,-2)  circle[radius=5pt];
			\fill[thick, black] (-2,0)  circle[radius=5pt];
			\fill[thick, black] (-2,2)  circle[radius=5pt];
			\node at (0.4,-0.6){$0$};
		\end{tikzpicture}
  	\caption{$X_3 = \CP^2 \# 3\overline{\CP^2}$}
	\end{subfigure}
\caption{The set $\mathcal{S}(\Delta)$ for the moment polytopes of the five toric del Pezzo surfaces.}
\label{fig:2dim}
\end{figure}
\noindent
\O bro and Paffenholz checked that all monotone polytopes in dimensions $\leqslant 9$ satisfy property $FS$. We therefore expect property $FS$ to hold for all monotone Delzant polytopes. The two-dimensional case is obvious by the classification of four-dimensional toric monotone symplectic manifolds, see Figure \ref{fig:2dim}. Let~$\Delta_0$ be the set of points $x\in \Delta$ such that $\min\{\ell_1(x),\ldots,\ell_n(x)\}$ is attained by exactly one $\ell_i(x)$. This is an open, dense subset of $\Delta$ which is subdivided into chambers $\Delta_i$ by the hyperplanes $\ell_i = \ell_j$, see Figure~\ref{fig:level1} in Section~\ref{sec:toric}. 

\begin{lemma}
\label{lem:polyenergy}
Let $(M,\omega)$ be a toric symplectic manifold whose moment polytope $\Delta$ satisfies property $FS$. Then 
	\begin{equation*}
		e_{\Delta}(x) = \min \{\ell_1(x),\dots,\ell_{k}(x)\}
	\end{equation*}
for all $x \in \Delta_0$.
\end{lemma}

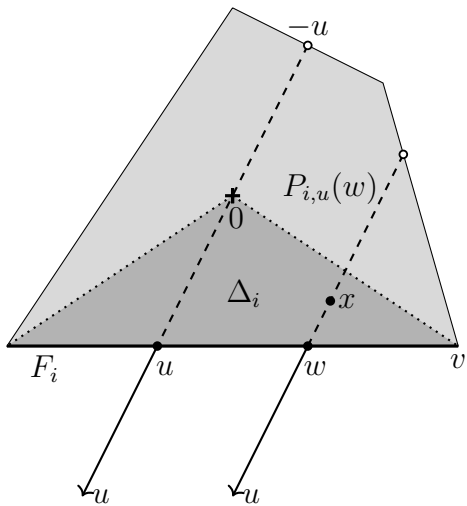
\begin{figure}[h]
	\begin{tikzpicture}[scale=1]
		\fill [black!15] (0,2.5)--(2,1.5)--(3,-2)--(-3,-2)--(0,2.5);
		\fill [black!30] (0,0)--(-3,-2)--(3,-2);
		\draw [thin,black] (0,2.5)--(2,1.5)--(3,-2)--(-3,-2)--(0,2.5);
		\draw[very thick] (-0.1,0)--(0.1,0)--(0.1,0);
		\draw[very thick] (0,0.1)--(0,-0.1)--(0,0.1);
		\draw[very thick] (-3,-2)--(3,-2);
		\draw[thick, dashed] (1,2)--(-1,-2);
		\draw[thick, dashed] (1,-2)--(2.27,0.55);
		\draw[thick, dotted] (0,0)--(-3,-2);
		\draw[thick, dotted] (0,0)--(3,-2);
		\draw[thick,->] (-1,-2)--(-2,-4);
		\draw[thick,->] (1,-2)--(0,-4);
		\fill[thick, black] (-1,-2)  circle[radius=1.8pt];
		\fill[thick, black] (1.3,-1.4)  circle[radius=1.8pt];
		\fill[thick, black] (1,2)  circle[radius=1.8pt];
		\fill[thick, black!0] (1,2)  circle[radius=1.1pt];
		\fill[thick, black] (1,-2)  circle[radius=1.8pt];
		\fill[thick, black] (2.27,0.55)  circle[radius=1.8pt];
		\fill[thick, black!0] (2.27,0.55)  circle[radius=1.1pt];
		\node at (-2.5,-2.3){$F_i$};
		\node at (-0.9,-2.3){$u$};
		\node at (-1.75,-4){$u$};
		\node at (0.25,-4){$u$};
		\node at (1,2.25){$-u$};
		\node at (1.1,-2.3){$w$};
		\node at (3,-2.2){$v$};
		\node at (1.52,-1.4){$x$};
		\node at (0.05,-0.3){$0$};
		\node at (1.3,0.1){$P_{i,u}(w)$};
		\node at (0.15,-1.3){$\Delta_i$};
		
	\end{tikzpicture}
	\caption{Construction of the probe $P_{i,u}(w)$.}
	\label{fig:probe2}
\end{figure}

\proof
The lower bound on the displacement energy follows from Proposition \ref{prop:lower}. For the upper bound, let $x \in \Delta_i$, which means that $\min\{\ell_1(x),\dots,\ell_k(x)\} = \ell_i(x)$. The set $\Delta_i$ is the cone $\{tx \,\vert\, t \in (0,1], \, x \in \mathring{F}_i\}$ over the interior~$\mathring{F}_i$ of the $i$-th facet~$F_i$ of~ $\Delta$. We are going to construct a probe with respect $F_i$ and apply Proposition \ref{prop:upper} for the upper bound. By the property $FS$, we can pick $u \in F_i \cap \mathcal{S}(\Delta)$. Since $u$ is integrally transverse to $F_i$ and $-u \in \Delta$, this yields a probe with the barycentre $0\in \Delta$ as its midpoint. Take the unique probe $P_{i,u}(w)$ parallel to $u$ which contains $x$, see Figure \ref{fig:probe2}. The point $x$ lies in the same half of $P_{i,u}(w)$ as $w$. This can be seen by noticing that the line segment between $-u$ and $v$ is contained in $\Delta$ by convexity of the moment polytope. 
\proofend

\section{Application I: Toric fibres}
\label{sec:toric}

An important class of examples for Lagrangian tori are moment fibres in toric symplectic manifolds. In this Section we use Theorem~\ref{t:dis} to give a criterion to exclude toric fibres from being real in terms of the function $e_{\Delta}$. We assume that $e_{\Delta}$ is given by the affine distance to the boundary of the moment polytope, see Assumption~\ref{ass:energy}. In Section~\ref{sec:displacement}, we proved that this assumption is reasonable in case the ambient manifold is monotone. In the present section, we do not assume monotonicity except for the proof of Theorem~\ref{t:toric}.\\ 

 Let $(M^{2n},\omega)$ be a compact toric symplectic manifold with moment map $\mu$ and moment polytope $\Delta$. For every point $x$ in the interior $\mathring{\Delta}$ of the moment polytope, the set $T_x = \mu^{-1}(x)$ is a Lagrangian torus in $M$ called toric fibre. Furthermore, the map 
\begin{eqnarray*}
	(H_1(T^n,\RR^n), 0) \cong (\RR^n,0)  &\rightarrow & \{\text{Lagrangian tori in }M \} \\
	a &\mapsto & T_{x+a} = \mu^{-1}(x+a),
\end{eqnarray*}
that is defined for all $a$ such that $x+a \in \mathring{\Delta}$, yields a versal deformation of $T_x$. Indeed, the components of $\mu$ give action coordinates on $\mu^{-1}(\mathring{\Delta})$ and thus $T_{x+a}$ and $T_{x+b}$ are related by a $C^1$-small Hamiltonian isotopy if and only if $a = b$. Varying $x$ in $\RR^n$ as above therefore yields an $n$-dimensional family of Hamiltonian isotopy classes of Lagrangian tori and hence a versal deformation of $T_x$.\\

\noindent
As a warm-up example and as an illustration to Theorem \ref{thm:main}, we consider the Clifford torus in products of $S^2$.

\begin{example}
\label{ex:sphereproduct}
{ \rm
Let $(S^2,\omega)$ be the unit $2$-sphere in $\RR^3$ equipped with the rescaled Euclidean area form $\omega = \frac{1}{2\pi} \text{area}$, for which $\int_{S^2} \omega = 2$. Let $H: S^2 \rightarrow \RR$ be the projection to the $z$-axis $H(p)=z$. Since the Hamiltonian flow of $H$ is $1$-periodic, it defines a toric structure on $S^2$ with moment polytope $[-1,1] \subset \RR$. The level sets of $T_c = H^{-1}(c)$ are circles of fixed height and have displacement energy 
\begin{equation*}
	e_{S^2}(T_c) = 
	\begin{cases} 
    	1 - \vert c \vert 	& \text{if } c \in [-1,1] \setminus \{0\} , \\
   		\infty				& \text{if } c = 0.
   \end{cases}
\end{equation*}
Recall from Example \ref{ex:equator} that the equator $T_0$ is real. In accordance with Theorem \ref{thm:main}, the displacement energy germ $S_{T_0}(c)=e_{S^2}(T_c)$ is invariant under $c\mapsto -c$. Consider the $n$-fold product of this example. The corresponding moment map $\mu$ is given as the $n$-fold product of the above Hamiltonian $H$. The moment polytope is the unit square $\Delta = [-1,1] \times \dots \times [-1,1] \subset \RR^n$. The level sets of $\mu$ are products of circles of fixed height. Their displacement energy is
\begin{equation*}
	e_{\times_n S^2}(T_{(c_1,\dots,c_n)}) =
	\begin{cases}
		\underset{1\leqslant i \leqslant n}{\min} \{ 1 - \vert c_i \vert \} 		& \text{if } (c_1,\dots,c_n) \in \Delta \setminus \{0\}; \\
		\infty														& \text{if } (c_1,\dots,c_n) = 0.
	\end{cases}
\end{equation*}
The \emph{Clifford torus} $T_0$ is real, and its displacement energy germ 
\begin{equation*}
S_{T_0} (c_1,\dots,c_n) = e_{\times_n S^2}(T_{(c_1,\dots,c_n)})
\end{equation*}
is invariant under $(c_1,\dots,c_n) \mapsto (-c_1,\dots,-c_n)$.
}
\end{example}

We will now turn to the class of toric symplectic manifolds for which the level sets of the function 
\begin{equation*}
e_{\Delta} \colon \Delta \rightarrow \RR \cup \{\infty\}, \quad x \mapsto e_M(T_x)
\end{equation*}
\noindent
look as in Figure~\ref{fig:level2} in Section~\ref{sec:intro}, namely like scalings of $\partial \Delta$. Let 
\begin{equation*}
\ell_i(x) = \kappa_i - \langle x , v_i \rangle 
\end{equation*}
be the functionals on $\RR^n$ which define $\Delta = \{\ell_i \geqslant 0,\, \forall i\}$, where the~$v_i$ are the primitive outward pointing normal vectors to the facets, see Subsection~\ref{ss:lbtoric}. The facets $F_i$ of $\Delta$ are given by the intersection of the moment polytope and the affine hyperplanes bounding the half-spaces, $F_i=\Delta \cap \{\ell_i=0\}$. For every $x \in \RR^n$, the value $\ell_i(x)$ is equal to the affine distance of $x$ to the corresponding facet $F_i$. See \cite{McD11} for details. 
\noindent
\begin{assumption} 
\label{ass:energy}
For all $x$ in an open dense subset of $\Delta$, the displacement energy of the toric fibre over $x\in \Delta$ is given by the affine distance of $x$ to the boundary $\partial\Delta$, i.e.\
\begin{equation*}
	e_{\Delta}(x) = \min \{ \ell_1(x),\dots,\ell_k(x) \}.
\end{equation*}
\end{assumption}
\noindent
If a variation of the Ewald conjecture holds, then this assumption is true for all monotone symplectic toric manifolds. See Section \ref{ssec:monotone} for details. \\

\begin{figure}
	\begin{tikzpicture}[scale=1.2]
		\fill [black!15] (0,0)--(1,3)--(4,0)--(0,0);
		\draw [thick,black,->] (1.5,0)--(1.5,-0.5);
		\draw [thick,black,->] (0.5,1.5)--(-1,2);
		\draw [thick, black] (1.5,0)--(1.6,0)--(1.6,-0.1)--(1.5,-0.1)--(1.5,0); 
		\draw [thick, black] (0.5,1.5)--(0.466,1.4)--(0.366,1.433)--(0.4,1.533)--(0.5,1.5); 
		\draw [thick] (1,3)--(0,0)--(4,0);
		\draw [thick,dashed,black] (0,0)--(2,2);
		\fill[thick, black] (1.98,0.8)  circle[radius=1.5pt];
		\fill[thick, black] (1.3,1.3)  circle[radius=1.5pt];
		\fill[thick, black] (1.2,2)  circle[radius=1.5pt];
		
		\node at (1.2,2.2){$x$};
		\node at (2.05,2.75){$I_x=\{j\}$};
		\node at (1.5,1.2){$y$};
		\node at (2.8,2.1){$I_y=\{i,j\}$};
		\node at (2.2,0.8){$z$};
		\node at (3.45,1.3){$I_z=\{i\}$};
		\node at (1.7,-0.4){$v_i$};
		\node at (2.8,-0.2){$F_i$};
		\node at (-0.7,1.7){$v_j$};
		\node at (0.6,2.5){$F_j$};
	\end{tikzpicture}
	\caption{The set $I_x$ for three different points.}
	\label{fig:level1}
\end{figure}
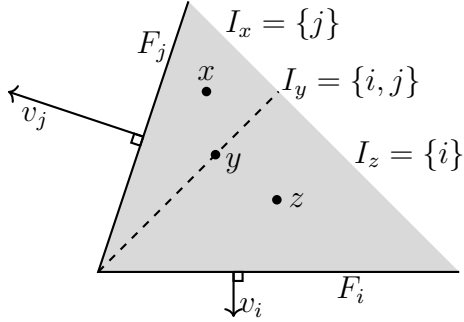

For any $x\in \mathring{\Delta}$ define the set $I_x$ of indices $i$ for which the minimal affine distance to $\partial \Delta$ is attained by the corresponding $\ell_i$, i.e. $i\in \{1,\dots,k\}$ belongs to $I_x$ if and only if $\ell_i(x) = \min \{ \ell_1(x), \dots \ell_k(x) \}$. Notice that if $I_x$ is not a singleton, then $x$ lies in a finite union of hyperplanes, see Figure \ref{fig:level1}.

\begin{proposition}
\label{prop:mingerm}
Let $x \in \Delta$. Under Assumption \ref{ass:energy}, the displacement energy germ of the corresponding toric fibre is given by
	\begin{equation*}
		S_{T_x}(a) = \min_{i \in I_x}\{\ell_i(x+a)\},
	\end{equation*}
for $a \in \RR^n$ in an open dense subset around $0$.
\end{proposition}
\proof By Assumption \ref{ass:energy}, we have
	\begin{equation*}
		S_{T_x}(a) = e(T_{x+a}) = \min_{1 \leq i \leq k} \{ \ell_i(x+a) \} = \min_{i \in I_x}\{\ell_i(x+a)\}.
	\end{equation*}
The last equality holds since $I_{x+a} \subseteq I_x$ for small enough $a$.
\proofend 

\begin{proposition}
\label{prop:symgerm}
Let $x \in \Delta$ be such that $T_x \subset M$ is a real Lagrangian. Under Assumption \ref{ass:energy} the moment polytope has to satisfy the following symmetry condition. For each $i\in I_x$ there is $j \in I_x$ such that $v_i = - v_j$. In particular, $I_x$ contains an even number of elements.
\end{proposition}
\proof
By Theorem \ref{thm:main}, if $T_x$ is real, then $S_{T_x}(a) = S_{T_x}(-a)$. By Proposition \ref{prop:mingerm}, this translates to 
	\begin{equation*}
		\min_{r \in I_x}\{\ell_r(x+a)\} = \min_{s \in I_x}\{\ell_s(x-a)\}
	\end{equation*}
for $a$ in an open neighbourhood of $0 \in \RR^n$. For every $i\in I_x$, there is an open set $U_i$ (which may not contain $0$) such that $\ell_i(x+a) = \min_{r \in I_x}\{\ell_r(x+a)\}$ for all $a \in U_i$. Hence, there is $j \in I_x$ such that, possibly after shrinking the subset $U_i$, we have
	\begin{equation*}
		\ell_i(x+a) = \ell_j(x-a), \quad \forall a \in U_i.
	\end{equation*}
Using $\ell_i(x+a)= \ell_i(x) - \langle v_i , a \rangle$, we deduce that $\langle v_i + v_j , a \rangle = 0$ for all $a \in U_i$ and hence $v_i = -v_j$.
\proofend

We are now in a position to prove Theorem~\ref{t:toric}.\\
\proofof{Theorem~\ref{t:toric}} Assume that $M$ is monotone and that its moment polytope satsifies property~$FS$. Lemma~\ref{lem:polyenergy} implies that Assumption~\ref{ass:energy} holds. Furthermore, since $M$ is monotone, we have $\ell_i(0)=1$ for all $i \in \{1,\ldots,k\}$ and thus $I_0 = \{1,\dots,k\}$. Hence the theorem follows from Proposition~\ref{prop:symgerm}. \proofend

\section{Application II: Chekanov tori} 
\label{sec:chekanov}

Chekanov tori were defined in~\cite{Che96} as the first examples of monotone Lagrangian tori in $\CC^n$ which are not symplectomorphic to a product torus. In this section, we recall an alternative construction given in~\cite{CheSch10}, see also~\cite{EliPol97}, and show that the Chekanov torus can be embedded into any toric monotone symplectic manifold. Under the property~$FS$, we compute its displacement energy germs and show that it is exotic and not real.

\subsection{Embedding Chekanov tori}
\label{ssec:chekemb}  Let $T^n$ act on $\CC^n$ by the standard Hamiltonian torus action generated by the moment map 
\begin{equation*}
	\nu \colon \CC^n \rightarrow \RR^n ,\quad (z_1,\ldots,z_n) \mapsto \pi (\vert z_1 \vert^2, \ldots , \vert z_n \vert^2 ) + (-1,\ldots,-1).
\end{equation*}
The image of $\nu$ is the positive quadrant in $\RR^n$ translated by the vector $(-1,\ldots,-1)$. By $\widehat{T}^{n-1}$ we will denote the linear subtorus 
\begin{equation}
	\label{eq:subtorus}
	\widehat{T}^{n-1} = \left\{ (e^{i\alpha_1},\ldots,e^{i\alpha_n}) \,\vert \, \alpha_1 + \ldots + \alpha_n = 0 \right\}\subset T^n,
\end{equation}
which has a natural Hamiltonian action on $\CC^n$. Now take a smooth embedded curve $\gamma(t) = r(t)e^{2\pi i\vartheta(t)}$ in $\CC$ which encloses area~$1$ and for which 
\begin{equation}
	\label{eq:gamma}
	0 < \vartheta(t) < \frac{1}{n} \quad \text{and} \quad 0 < r(t) < \sqrt{\frac{n}{\pi}} + \delta,
\end{equation}
for a small $\delta > 0$. From $\gamma$ construct the curve $\Gamma(t) = \frac{1}{\sqrt{n}}(\gamma(t),\ldots,\gamma(t))$ lying in the diagonal plane in $\CC^n$.

\begin{figure}
		\centering
		\begin{tikzpicture}[scale=3]
			\fill[black!15] (-1,0.5)--(-1,-1)--(0.5,-1)--(0.5,0.5)--(-1,0.5);
			\fill[black!35] (0.2,-0.1) arc (-45:135:0.212);
			\fill[black!35] (0.2,-0.1)--(-0.7,-1)--(-1,-1)--(-1,-0.7)--(-0.1,0.2)--(0.2,-0.1);
			\draw[thin,->] (0,-1.3)--(0,0.8);
			\draw[thin,->] (-1.3,0)--(0.8,0);
			\draw[thick] (-1,0.5)--(-1,-1)--(0.5,-1);
			\draw[thick,dotted] (-1,-1)--(0.5,0.5);
			\draw[thick] (-0.9,-0.9)--(0.1,0.1);
			\node at (0.14,0.8){$x_2$};
			\node at (0.8,0.10){$x_1$};
			\node at (-1,-1.15){$(-1,-1)$};
			\node at (-0.3,-0.5){$\nu(\Theta^n)$};
			\node at (-0.2,0.2){$U$};
			\fill[thick, black] (-0.9, -0.9)  circle[radius=0.75pt];
			\fill[thick, black] (0.1, 0.1)  circle[radius=0.75pt];
		\end{tikzpicture}	
\caption{The image of $\Theta^n \subset \CC^n$ under $\nu$.}
\label{fig:thetamoment}
\end{figure}
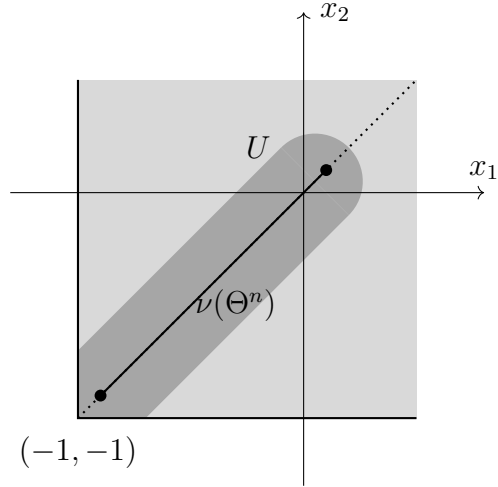

\begin{definition}
The \emph{Chekanov torus} $\Theta^n$ in $\CC^n$ is the torus swept out by $\Gamma$ under the action of $\widehat{T}^{n-1}$,
	\begin{equation*}
		\Theta^n = \left\{\left. \frac{1}{\sqrt{n}} \left(e^{i\alpha_1} \gamma(t),\ldots, e^{i\alpha_n}\gamma(t)   \right) \in \CC^n \right\vert \alpha_1 + \ldots + \alpha_n = 0  \right\}.
	\end{equation*}
\end{definition}
\noindent
The Chekanov torus is embedded, Lagrangian and monotone. Notice that $\nu(\Theta^n)$ is contained in the diagonal line, and by the choice of~$\gamma$ in~$(\ref{eq:gamma})$ every component satisfies 
\begin{equation}
	\label{eq:nui}
	\varepsilon - 1 < \nu_i(\Theta^n) < \varepsilon
\end{equation} 
for a small $\varepsilon > 0$, see Figure~\ref{fig:thetamoment}.

\begin{remark}
{\rm
The Chekanov torus $\Theta^n \subset \CC^n$ and is not real. In fact, by the Smith inequality~$(\ref{eq:smith2})$, tori in $\CC^n$ cannot be realized as the fixed point set of a \emph{smooth} involution.
}
\end{remark}

Let $M^{2n}$ be a toric monotone symplectic manifold with moment map~$\mu$. We show that $\Theta^n$ can be embedded into $M$. Pick a vertex $v$ of its moment polytope $\Delta = \mu(M) = \{\ell_i \geqslant 0\}$. Since $\Delta$ is a Delzant polytope, we can assume (up to applying a transformation in $\GL(n,\ZZ)$) that the facets meeting at $v$ are parallel to the coordinate hyperplanes. By monotonicity, these hyperplanes lie at affine distance $1$ to the origin and hence $v = (-1,\ldots,-1)$. In other words, we assume that the $n$ first functionals defining $\Delta$ satisfy
\begin{equation}
	\label{eq:ellnormal}
	 \ell_1(x) = 1 + x_1, \, \ldots \, , \ell_n(x) = 1 +x_n.
\end{equation}
Equivalently, the moment polytope $\Delta$ near $v$ has the same structure as $\nu(\CC^n)$ near $\nu(0)$. This can be used to construct an embedding of $\Theta^n$ into $M$. By convexity of $\Delta$, the line segment between the origin and $v$ is contained in $\Delta$. By~$(\ref{eq:nui})$ we can thus choose a neighbourhood $U \subset \nu(\CC^n)$ of the segment $\nu(\Theta^n)$ which fits into $\Delta$, see again Figure~\ref{fig:thetamoment}. This yields a $T^n$-equivariant symplectic embedding of the neighbourhood $\nu^{-1}(U)$ of $\Theta^n$ into $M$. Denote the so obtained Chekanov torus by $\Theta^n_M$. By equivariance of the embedding, it is invariant under the $\widehat{T}^{n-1}$-action on $M$ induced by $\mu$.

\begin{proposition}
\label{prop:Chekembed}
Let $M$ be a toric monotone symplectic manifold. Then the Chekanov torus embeds into $M$ to yield a monotone Lagrangian torus $\Theta^n_M \subset M$.
\end{proposition}

\proof We prove that~$\Theta^n_M$ is monotone. This means that the Maslov index and the area class are proportional on disks with boundary on $\Theta^n_M$, i.e.\ that there is a~$C > 0$ such that
	\begin{equation*}
		\maslov(D) = C \int_D \omega, \quad \forall D \in \pi_2(M,\Theta^n_M).
	\end{equation*}
The homotopy long exact sequence yields
	\begin{equation*}
		0 \rightarrow \pi_2(M) \rightarrow \pi_2(M,\Theta^n_M) \rightarrow \pi_1(\Theta^n_M) \rightarrow 0.
	\end{equation*}
As a basis for~$\pi_1(\Theta^n_M)$ we choose~$[\Gamma]$ and the orbits of the~$\widehat{T}^{n-1}$-action. The Maslov index and the area class vanish on the latter and $\maslov([\Gamma])=2 \int_{[\Gamma]}\omega = 2$. On spheres, the Maslov index is equal to twice the first Chern class of the ambient manifold. Recall that~$M$ is itself monotone with $c_1 = [\omega]$ and thus we obtain $\maslov(D) = 2c_1(D) = 2\int_{D}\omega$ for all $D \in \pi_2(M)$. This proves that~$\Theta^n_M$ is monotone with $C=2$. 
\proofend

\begin{remark}
{\rm
In general, the tori $\Theta_M$ may depend on the choice of the vertex $v$. However, in the cases of $\times_n S^2$ and $\CP^n$ all vertices of the corresponding moment polytopes are interchangeable by an element of $\GL(n,\ZZ)$ and hence we obtain a unique torus $\Theta_M$ up to symplectomorphism. 
}
\end{remark}

\subsection{Versal deformations}
Assume that $M$ has property $FS$. For readability, we will write $\stackrel{\otimes}{=}$ for equalities that hold on an open dense subset of a neighbourhood of the origin of a vector space. By monotonicty, we can assume~$\ell_i(0)=1$ for all~$i$ and hence, by Proposition~\ref{prop:mingerm},
\begin{equation*}
	S_{T_0}(a) = e_{\Delta}(a) \stackrel{\otimes}{=} \min\{\ell_1(a),\ldots,\ell_n(a)\}.
\end{equation*}
In particular, the displacement energy germ of the central fibre~$T_0$ is determined by the moment polytope. The displacement energy germ of the corresponding Chekanov torus $\Theta_M$ is closely related to the one of~$T_0$.

\begin{lemma}
\label{lem:thetaversal}
Let $M$ be a toric monotone symplectic manifold satisfying property $FS$. Then the displacement energy germ of the Chekanov torus $\Theta^n_M$ is given by
\begin{equation}
	\label{eq:germ1}
	S_{\Theta^n_M} \stackrel{\otimes}{=} S_{T_0} \circ \phi.
\end{equation}
Here $\phi \colon \RR^n \rightarrow \RR^n$ is the piece-wise linear homeomorphism defined by~$(\ref{eq:phi1})$ and~$(\ref{eq:phi2})$, which does not depend on $M$.
\end{lemma}

\proof We will closely follow the ideas used in~\cite{CheSch10} to compute $S_{\Theta_{S^2\times S^2}}$. Since there is no risk of confusion here, we denote the Chekanov torus by $\Theta = \Theta^n_M$. Let $\mu \colon M \rightarrow \RR^n$ be the moment map for which $\Delta$ has the form~$(\ref{eq:ellnormal})$. Notice that the subtorus $\widehat{T}^{n-1}$ defined by equation~$(\ref{eq:subtorus})$ has a natural Hamiltonian action on $M$ via the inclusion $\widehat{T}^{n-1} \subset T^n$ and that $\Theta$ is invariant under this torus action. The moment map $\widehat{\mu} \colon M \rightarrow \RR^{n-1}$ corresponding to the $\widehat{T}^{n-1}$-action is given by 
\begin{equation}
	\label{eq:muhat}
	\widehat{\mu} = (\mu_1 - \mu_n , \ldots , \mu_{n-1} - \mu_n).
\end{equation}
As a basis of $H_1(\Theta,\ZZ)$, we choose the class~$[\Gamma]$ of the curve lying in the diagonal and the classes $[\tau_1],\ldots,[\tau_{n-1}]$ of the orbits of the $\widehat{T}^{n-1}$-action. The latter can also be seen as the closed orbits of the Hamiltonians~$\mu_i - \mu_n$. By the equivariant Weinstein neighbourhood theorem, we can choose a versal deformation of $\Theta$ which preserves the $\widehat{T}^{n-1}$-orbit structure. Let $t_1,\ldots,t_{n-1}$ and $s$ be the deformation parameters corresponding to the classes $[\tau_1],\ldots,[\tau_{n-1}]$ and $[\Gamma]$. For convenience we denote $\mathbf{t}=(t_1,\ldots,t_n) \in \RR^{n-1}$. Since $\widehat{T}^{n-1}$-orbits are preserved, we find that the Lagrangian neighbour $\Theta_{\mathbf{t},s}$ of $\Theta$ maps to a line segment $\mu(\Theta_{\mathbf{t},s})$ parallel to $\mu(\Theta)$. Furthermore, by equation~$(\ref{eq:muhat})$, the line segment $\mu(\Theta_{\mathbf{t},s})$ is contained in the line
\begin{equation}
	\label{eq:lt}
	L_{\mathbf{t}}= \left\{(x_1,\ldots,x_n) \in \RR^n \, \vert \, x_1 - x_n = t_1 ,\ldots, x_{n-1} - x_n = t_{n-1}\right\}.
\end{equation}

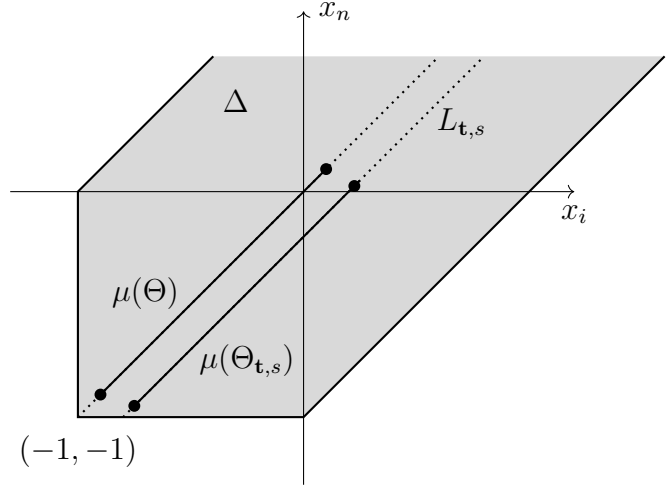
\begin{figure}
		\centering
		\begin{tikzpicture}[scale=3]
			\fill[black!15] (-1,-1)--(0,-1)--(1.6,0.6)--(0.6,0.6)--(-0.4,0.6)--(-1,0)--(-1,-1);
			\draw[thin,->] (0,-1.3)--(0,0.8);
			\draw[thin,->] (-1.3,0)--(1.2,0);
			\draw[thick] (-0.4,0.6)--(-1,0)--(-1,-1)--(0,-1)--(1.6,0.6);
			\draw[thick,dotted] (-1,-1)--(0.6,0.6);
			\draw[thick,dotted] (-0.8,-1)--(0.8,0.6);
			\draw[thick] (-0.9,-0.9)--(0.1,0.1);
			\draw[thick] (-0.75,-0.95)--(0.25, 0.05);
			\node at (0.14,0.8){$x_n$};
			\node at (1.2,-0.10){$x_i$};
			\node at (-1,-1.15){$(-1,-1)$};
			\node at (-0.7,-0.45){$\mu(\Theta)$};
			\node at (-0.25,-0.75){$\mu(\Theta_{\mathbf{t},s})$};
			\node at (0.7,0.3){$L_{\mathbf{t},s}$};
			\node at (-0.3,0.4){$\Delta$};
			\fill[thick, black] (-0.9, -0.9)  circle[radius=0.75pt];
			\fill[thick, black] (0.1, 0.1)  circle[radius=0.75pt];
			\fill[thick, black] (-0.75, -0.95)  circle[radius=0.75pt];
			\fill[thick, black] (0.225, 0.025)  circle[radius=0.75pt];
		\end{tikzpicture}	
\caption{Versal deformation of $\Theta \subset M$.}
\label{fig:thetamoment2}
\end{figure}
\noindent 
See Figure~\ref{fig:thetamoment2}. We prove that whenever $t_i \neq 0$ for all $1\leqslant i \leqslant n-1$, the versal deformation $\Theta_{\mathbf{t},s}$ of $\Theta$ is Hamiltonian isotopic to a toric fibre $T_{x} = \mu^{-1}(x)$ for a suitable $x = (x_1,\ldots,x_n)$. Since the displacement energy is preserved under Hamiltonian isotopies, property~$FS$ and Lemma~\ref{lem:polyenergy} yield the displacement energy germ of $\Theta$. Notice that if~$t_i \neq 0$ for all~$i$, then~$\widehat{T}^{n-1}$ acts freely on the set 
\begin{equation*}
Z_{\mathbf{t}} = \mu^{-1}(L_{\mathbf{t}} \cap \Delta \setminus \{y\} ) = \widehat{\mu}^{-1}(\mathbf{t}) \setminus \mu^{-1}(y)
\end{equation*}
since $L_{\mathbf{t}}$ hits the boundary $\partial \Delta$ in a codimension one facet\footnote{Since the opposite point $y$ in the intersection $L_{\mathbf{t}} \cap \Delta$ might lie in a face with higher codimension, we remove it in order for the action to be free.}. Hence, we can perform symplectic reduction by $\widehat{T}^{n-1}$ on $Z_{\mathbf{t}}$ 
\begin{center}
	\begin{tikzcd}
		\Theta_{\mathbf{t},s} \subset Z_{\mathbf{t}} \arrow[hook]{r}{i} \arrow[two heads]{d}{p}
		& (M,\omega) \\
		c_{\mathbf{t},s} \subset (M_{\mathbf{t}},\omega_{\mathbf{t}}).
	\end{tikzcd}
\end{center}
The symplectic quotient $(M_{\mathbf{t}},\omega_{\mathbf{t}})$ is symplectomorphic to a disk of radius equal to the affine length of $L_{\mathbf{t}} \cap \Delta$. Indeed, since there is a Hamiltonian $T^n$-action on $Z_{\mathbf{t}}$ the reduced space has an induced Hamiltonian $S^1$-action with moment polytope $L_{\mathbf{t}} \cap \Delta \setminus \{y\}$. Since $\Theta_{\mathbf{t},s}$ is $\widehat{T}^{n-1}$-invariant, it projects to a circle $c_{\mathbf{t},s} = p (\Theta_{\mathbf{t},s})$. We claim that this circle encloses symplectic area $ 1+ s$. Since the $\widehat{T}^{n-1}$-orbits $\tau_1,\ldots,\tau_{n-1}$ are divided out by the above symplectic reduction, the circle $c_{\mathbf{t},s}$ corresponds to the class $[\Gamma]$ in $\Theta_{\mathbf{t},s}$. The latter class bounds a disk of area $1+s$ in $M$ since $s$ is the deformation parameter of $[\Gamma]$. By symplectic reduction we have $p^*\omega_{\mathbf{t}} = i^*\omega$ and hence $c_{\mathbf{t},s}$ encloses area $1+s$. It is thus Hamiltonian isotopic to the concentric circle $S^1(1+s)$ in $M_{\mathbf{t}}$ which bounds the same area. The pre-image $p^{-1}(S^1(1+s))$ is a toric fibre $T_x$ and thus the Hamiltonian isotopy in the quotient can be lifted to $M$ to yield a Hamiltonian isotopy between $\Theta_{\mathbf{t},s}$ and $T_x$. 

Now, let $\phi \colon \RR^n \rightarrow \RR^n$ be the map that takes $(\mathbf{t},s)$ to $x$ such that $\Theta_{\mathbf{t},s}$ and $T_x$ are Hamiltonian isotopic. Note that this defines $\phi$ only on an open dense subset of a neighbourhood of $0 \in \RR^n$ on which we have
\begin{equation*}
	\label{eq:germ3}
	S_{\Theta}(\mathbf{t},s) 
	= e_M(\Theta_{\mathbf{t},s})
	\stackrel{\otimes}{=} e_M(T_{\phi(\mathbf{t},s)}) 
	= e_{\Delta}(\phi(\mathbf{t},s))
	\stackrel{\otimes}{=} S_{T_0}(\phi(\mathbf{t},s)).
\end{equation*}
We now determine the map $\phi$. Let $(\mathbf{t},s) \in \RR^n$ be such that $\phi$ is defined. The point $x = \phi(\mathbf{t},s)$ lies on $L_{\mathbf{t}}$ and hence 
\begin{equation*}
	\label{eq:xt1}
	t_1 = x_1 - x_n,\ldots, t_{n-1} = x_{n-1} - x_n.
\end{equation*}
Let $z\in \partial \Delta$ be the point close to $(-1,\ldots,-1)$ in which $L_{\mathbf{t}}$ intersects the boundary of $\Delta$. The area enclosed by $S^1(1+s) \subset M_{\mathbf{t}}$ is equal to the affine length of the line segment~$[z,x]$, which in turn is equal to $1+ \min\{x_1,\ldots,x_n\}$ and hence 
\begin{equation}
	\label{eq:xt2}
	s = \min\{x_1,\ldots,x_n\}.
\end{equation}
The map $\phi$ we are looking for is thus given as the inverse of 
\begin{equation}
	\label{eq:xt3}
	\begin{pmatrix}
		x_1 \\ 
		\vdots \\
		x_{n-1} \\ 
		x_n
	\end{pmatrix}
	\mapsto
	\begin{pmatrix}
		x_1 - x_n \\ 
		\vdots \\
		x_{n-1} -x_n \\ 
		\min\{x_1,\ldots,x_n\}
	\end{pmatrix}.
\end{equation}
There is a unique extension to a piece-wise linear homeomorphism on all of $\RR^n$. By distinguishing cases we obtain 
\begin{equation}
	\label{eq:phi1}
	\phi(\mathbf{t},s)=
	\begin{pmatrix}
		 	s + t_1 \\
		 	\vdots \\
		 	s + t_{n-1} \\
		 	s
	\end{pmatrix},
\end{equation}
whenever all $t_i \geqslant 0$ and
\begin{equation}
	\label{eq:phi2}
	\phi(\mathbf{t},s)=
	\begin{pmatrix}
		 	s + t_1 - t_i \\
		 	\vdots \\
		 	s + t_{n-1} - t_i \\
		 	s - t_i
	\end{pmatrix},
\end{equation}
if $t_i < 0$ and $t_i$ is minimal among all $t_j$.
\proofend

Instead of working directly with the displacement energy germ $S_L$ of a Lagrangian $L$, it is often useful to look at its level sets $S_L^{-1}(c)$ for some $c>0$. In particular, if $L$ is real, then these level sets are centrally symmetric, by Theorem~\ref{thm:main}. In the case of $T_0$, the level sets are rescalings of $\Delta$,
\begin{equation}
	\label{eq:germlevel2}
	S_{T_0}^{-1}(c) \stackrel{\otimes}{=} \lambda \Delta, \quad \lambda > 0.
\end{equation}
Here we mean that both sets agree when intersected with a set which is open and dense in the neighbourhood of the origin. Since $S_{\Theta^n_M} \stackrel{\otimes}{=} S_{T_0} \circ \phi$, we obtain 
\begin{equation}
	\label{eq:germlevel}
	S_{\Theta^n_M}^{-1}(c) 
	\stackrel{\otimes}{=} \phi^{-1}(S_{T_0}^{-1}(c))
	\stackrel{\otimes}{=} \lambda \phi^{-1}(\Delta).
\end{equation}
This allows us to understand the versal deformation of $\Theta^n_M$ by applying $\phi^{-1}$ to the moment polytope $\Delta$. The inverse of $\phi$ is given by equation~$(\ref{eq:xt3})$.

We will now prove that one can pick a suitable vertex $v$ for which the embedding of the Chekanov torus constructed in Subsection \ref{ssec:chekemb} yields an exotic Lagrangian torus in $M$, i.e.\ a torus which is not symplectomorphic to a toric fibre. For this, let $F_0$ be a facet of the moment polytope $\Delta$ which contains the maximal number of integral points among all facets of $\Delta$, let $v$ to be any vertex contained in $F_0$ and let $\Theta_M^n$ be the Chekanov torus embedded with respect to $v$. A priori, $\Theta_M^n$ can only be symplectomorphic to the central toric fibre, since all other fibres are not monotone. By \eqref{eq:germlevel2} and \eqref{eq:germlevel}, it suffices to show that the polytopes~$\Delta$ and $\phi^{-1}(\Delta)$ are not $\GL(n,\ZZ)$-equivalent in order to show that~$T_0$ and ~$\Theta_M^n$ are not symplectomorphic. Note that the maximal number of lattice points in a facet is a $\GL(n,\ZZ)$-invariant of polytopes and thus it suffices to show that this invariant strictly increases when we apply~$\phi^{-1}$ with respect to~$v$. Assume that $\Delta$ is given in the normal form~\eqref{eq:ellnormal} with respect to $v$ and hence the minimum $\min\{x_1,\ldots,x_n\}$ is constant and equal to~$-1$ on all facets containing~$v=(-1,\ldots,-1)$. Therefore $\phi^{-1}$ maps all facets containing~$v$ (in particular~$F_0$) to the same facet of~$\phi^{-1}(\Delta)$, which therefore contains strictly more integral points than any facet in $\Delta$. We have shown

\begin{proposition}
\label{prop:exotic}
Let $M$ be a toric monotone symplectic manifold satisfying property $FS$. Then $M$ contains an exotic copy of the Chekanov torus.
\end{proposition}

\begin{remark}
{\rm
The following example shows that the right choice of the vertex $v$ is crucial for the obtained Chekanov torus to be distinguishable from the central fibre by versal deformations. The polytope in $\RR^2$ defined by the functionals
\begin{equation*}
	1 + x_1, \; 1 \pm x_2, \; 1-x_1+x_2
\end{equation*}
is the moment polytope of the one-fold blow-up $X_1$ of $\CP^2$. The level sets of $S_{\Theta^2_{X_1}}$ when $\Theta^2$ is embedded with respect to the vertex $(-1,-1)$ are rescalings of the polytope defined by
\begin{equation*}
	1 - t, \; 1 \pm s, \; 1 + t - s.
\end{equation*}
Since these two polytopes are related by an element in $\GL(2,\ZZ)$, versal deformations cannot distinguish between~$T^2_0 \subset X_1$ and~$\Theta^2_{X_1}$.
}
\end{remark}

\begin{figure}
	\begin{subfigure}{0.45\textwidth}
  		\centering
		\begin{tikzpicture}[scale=0.6]
			\draw[thick,->] (0,-2.5)--(0,2.5);
			\draw[thick,->] (-2.5,0)--(2.5,0);
			\draw[thin] (-4,-2)--(4,-2)--(0,2)--(-4,-2);
			\draw[thin] (-3.5,-1.75)--(3.5,-1.75)--(0,1.75)--(-3.5,-1.75);
			\draw[thin] (-3,-1.5)--(3,-1.5)--(0,1.5)--(-3,-1.5);
			\draw[thin] (-2.5,-1.25)--(2.5,-1.25)--(0,1.25)--(-2.5,-1.25);
			\draw[thin] (-2,-1)--(2,-1)--(0,1)--(-2,-1);
			\draw[thin] (-1.5,-0.75)--(1.5,-0.75)--(0,0.75)--(-1.5,-0.75);
			\draw[thin] (-1,-0.5)--(1,-0.5)--(0,0.5)--(-1,-0.5);
			\draw[thin] (-0.5,-0.25)--(0.5,-0.25)--(0,0.25)--(-0.5,-0.25);
			\node at (0.35,2.5){$s$};
			\node at (2.5,0.375){$t$};
		\end{tikzpicture}
  	\caption{$M = S^2 \times S^2$}
	\end{subfigure}
	\begin{subfigure}{0.45\textwidth}
		\centering
		\begin{tikzpicture}[scale=0.6]
			\draw[thick,->] (0,-2.5)--(0,2.5);
			\draw[thick,->] (-2.5,0)--(2.5,0);
			\draw[thin] (-4,-2)--(6,-2)--(0,1)--(-2,0)--(-4,-2);
			\draw[thin] (-3.5,-1.75)--(5.25,-1.75)--(0,0.875)--(-1.75,0)--(-3.5,-1.75);
			\draw[thin] (-3,-1.5)--(4.5,-1.5)--(0,0.75)--(-1.5,0)--(-3,-1.5);
			\draw[thin] (-2.5,-1.25)--(3.75,-1.25)--(0,0.625)--(-1.25,0)--(-2.5,-1.25);
			\draw[thin] (-2,-1)--(3,-1)--(0,0.5)--(-1,0)--(-2,-1);
			\draw[thin] (-1.5,-0.75)--(2.25,-0.75)--(0,0.375)--(-0.75,0)--(-1.5,-0.75);
			\draw[thin] (-1,-0.5)--(1.5,-0.5)--(0,0.25)--(-0.5,0)--(-1,-0.5);
			\draw[thin] (-0.5,-0.25)--(0.75,-0.25)--(0,0.125)--(-0.25,0)--(-0.5,-0.25);
			\node at (0.35,2.5){$s$};
			\node at (2.5,0.375){$t$};
		\end{tikzpicture}
	\caption{$M = X_1$}
	\end{subfigure}
	\begin{subfigure}{0.45\textwidth}
  		\centering
		\begin{tikzpicture}[scale=0.6]
			\draw[thick,->] (0,-2.5)--(0,2.5);
			\draw[thick,->] (-2.5,0)--(2.5,0);
			\draw[thin] (-2,-2)--(2,-2)--(2,0)--(0,2)--(-2,0)--(-2,-2);
			\draw[thin] (-1.75,-1.75)--(1.75,-1.75)--(1.75,0)--(0,1.75)--(-1.75,0)--(-1.75,-1.75);
			\draw[thin] (-1.5,-1.5)--(1.5,-1.5)--(1.5,0)--(0,1.5)--(-1.5,0)--(-1.5,-1.5);
			\draw[thin] (-1.25,-1.25)--(1.25,-1.25)--(1.25,0)--(0,1.25)--(-1.25,0)--(-1.25,-1.25);
			\draw[thin] (-1,-1)--(1,-1)--(1,0)--(0,1)--(-1,0)--(-1,-1);
			\draw[thin] (-0.75,-0.75)--(0.75,-0.75)--(0.75,0)--(0,0.75)--(-0.75,0)--(-0.75,-0.75);
			\draw[thin] (-0.5,-0.5)--(0.5,-0.5)--(0.5,0)--(0,0.5)--(-0.5,0)--(-0.5,-0.5);
			\draw[thin] (-0.25,-0.25)--(0.25,-0.25)--(0.25,0)--(0,0.25)--(-0.25,0)--(-0.25,-0.25);
			\node at (0.35,2.5){$s$};
			\node at (2.5,0.375){$t$};
		\end{tikzpicture}
  	\caption{$M = X_3$}
	\end{subfigure}
\caption{Level sets of the function $S_{\Theta^n_M}$.}
\label{fig:level3}
\end{figure}
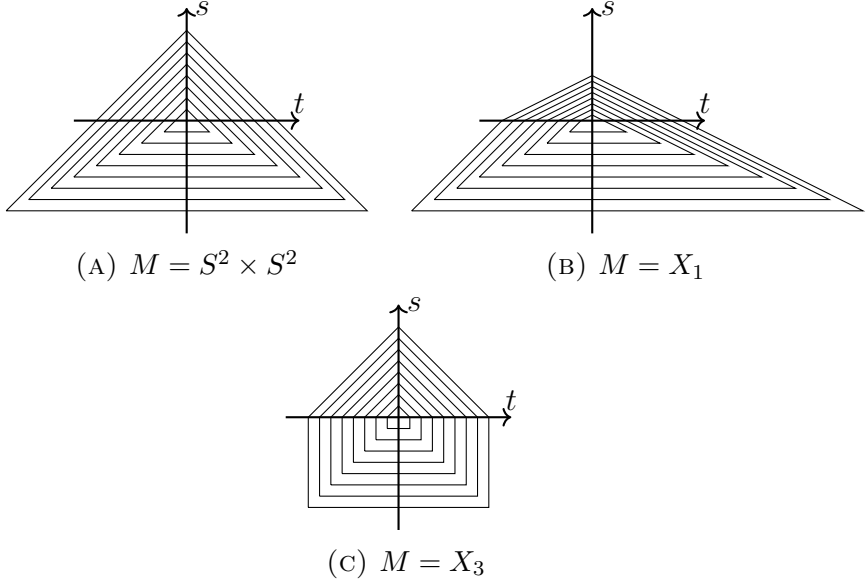

\subsection{Chekanov tori are not real} As a warm-up, let $M \in \{S^2\times S^2, \CP^2, X_1,X_2,X_3\}$ be one of the five toric monotone symplectic manifolds in dimension $4$. See Figure \ref{fig:2dim} in Section \ref{sec:displacement} for their moment polytopes. Then $\Theta^2_M$ is not real. The existence of real Lagrangian tori in $M=\CP^2$ and $M=X_2$ is excluded by the Smith inequality, see~($\ref{eq:smith2}$) in Section~\ref{sec:intro}. Applying $\phi^{-1}$ to the moment polytopes of the remaining three cases shows that the corresponding Chekanov tori are not real either, since the level sets of their displacement energy germs are not centrally symmetric, see Figure~\ref{fig:level3}. This can be generalized to all $\Theta^n_M$.

\begin{theorem}
\label{thm:chek}
Let $M$ be a toric monotone symplectic manifold satisfying property $FS$. Then the Chekanov torus $\Theta^n_M$ is not real.
\end{theorem}

\begin{figure}
	\begin{subfigure}{0.4\textwidth}
		\centering
		\begin{tikzpicture}[scale=1.5]
			\fill[thick, black!20] (-1, -1)  circle[radius=8.45pt];
			\fill[thick, black!20] (1, 1)  circle[radius=8.45pt];
			\fill[thin,black!0] (-1,-1)--(-0.7,-1)--(-0.7,-1.3)--(-1.3,-1.3)--(-1.3,-0.7)--(-1,-0.7);
			\fill[thin,black!0] (1,1)--(1,1.3)--(1.3,1.3)--(1.3,1);
			
			\draw[step=1.0,black!30,thin] (-1.3,-1.5) grid (1.3,1.5);
			\node at (0.25,1.5){$x_n$};
			\node at (1.3,-0.2){$x_1$};
			\node at (-0.4,-0.7){$U \cap \Delta$};
			\node at (0.75,0.55){$\phi(-V) \cap \Delta$};
			
			\draw[very thick] (-1,-0.7)--(-1,-1)--(-0.7,-1);
			\fill[thick, black] (-1, -1)  circle[radius=0.75pt];
			
			\draw[very thick] (1,1.3)--(1,1)--(1.3,1);
			\fill[thick, black] (1, 1)  circle[radius=0.75pt];
			
			\draw[thin,->] (0,-1.5)--(0,1.5);
			\draw[thin,->] (-1.3,0)--(1.3,0);

		\end{tikzpicture}	
		\caption{$\Delta$}
	\end{subfigure}
	\begin{subfigure}{0.1\textwidth}
		\centering
		\begin{tikzpicture}[scale=1.5]			
			\draw[thin,->] (-0.3,0)--(0.3,0);
			\node at (0,0.2){$\phi^{-1}$};
		\end{tikzpicture}	
	\end{subfigure}
	\begin{subfigure}{0.4\textwidth}
		\centering
		\begin{tikzpicture}[scale=1.5]
			\fill[thick, black!20] (0, -1)  circle[radius=8.45pt];
			\fill[thick, black!20] (0, 1)  circle[radius=8.45pt];
			\fill[thin,black!0] (-1,-1)--(1,-1)--(1,-1.3)--(-1,-1.3);
			\fill[thin,black!0] (-1,1)--(1,1)--(1,1.3)--(-1,1.3);
			
			\draw[step=1.0,black!30,thin] (-1.3,-1.5) grid (1.3,1.5);
			\draw[thin,->] (0,-1.5)--(0,1.5);
			\draw[thin,->] (-1.3,0)--(1.3,0);
			\node at (0.25,1.5){$x_n$};
			\node at (1.3,-0.2){$x_1$};
			\node at (0.85,-0.6){$V \cap \phi^{-1}(\Delta)$};
			\node at (0.85,0.5){$-V \cap \phi^{-1}(\Delta)$};
			
			\draw[very thick] (-0.3,-1)--(0.3,-1);
			\fill[thick, black] (0, -1)  circle[radius=0.75pt];
			
			\draw[very thick] (-0.3,1)--(0.3,1);
			\fill[thick, black] (0, 1)  circle[radius=0.75pt];
			
			\draw[thin,->] (0,-1.5)--(0,1.5);
			\draw[thin,->] (-1.3,0)--(1.3,0);
		\end{tikzpicture}	
	\caption{$\phi^{-1}(\Delta)$}
	\end{subfigure}
\caption{Idea of the proof of Theorem~\ref{thm:chek}, the set $\phi(-V) \cap \Delta$ is not convex.}
\label{fig:nreal}
\end{figure}
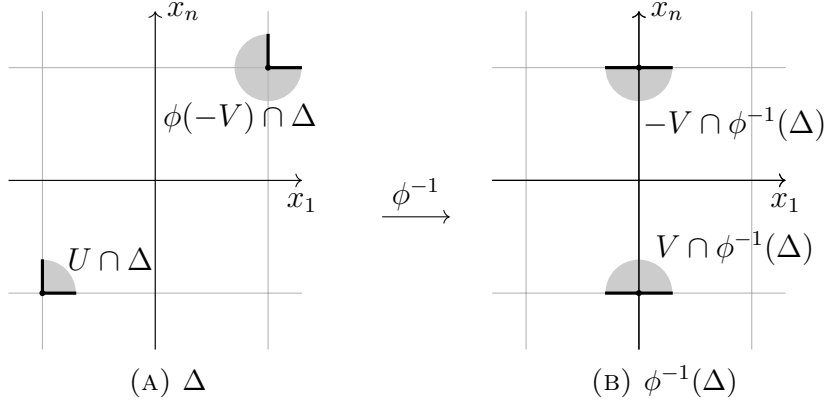

\proof Again, we suppose that $\Delta$ is in the form~$(\ref{eq:ellnormal})$ with distinguished vertex~$v=(-1,\ldots,-1)$. In order to understand the versal deformation of $\Theta^n_M$, we apply $\phi^{-1}$ to the moment polytope as in~$(\ref{eq:germlevel})$. The vertex~$v$ is mapped to $-e_n$ and all facets surrounding it to the hypersurface~$\{s=-1\}$. Hence, if $U \subset \RR^n$ is a neighbourhood of $v$, then there is a neighbourhood $V \subset \RR^n$ of $-e_n$ such that
\begin{equation*}
	\phi^{-1}(U \cap \partial \Delta) = V \cap \{s=-1\} \subset \partial \phi^{-1}(\Delta).
\end{equation*}
Now suppose that $\Theta^n_M$ is real and hence, by~$(\ref{eq:germlevel})$ that $\phi^{-1}(\Delta)$ is centrally symmetric. This implies that 
\begin{equation*}
	(-V) \cap \{s=1\} \subset \partial \phi^{-1}(\Delta)
\end{equation*}
Since~$-V$ is a neighbourhood of~$e_n$, points of the form~$e_n + re_i$ belong to~$(-V) \cap \{s=1\}$ and hence to $\phi^{-1}(\Delta)$ for small~$r>0$ and~$i\neq n$. This implies that $\phi(e_n + re_i) \in \Delta$. Observe that $\phi(e_n + re_i) = (1,\ldots,1) + re_i $ by equation~$(\ref{eq:phi1})$. Since $(1,\ldots,1)$ is integral, it does not belong to the interior of $\Delta$ and hence $\phi(e_n + re_i) \in \Delta$ contradicts the convexity of the moment polytope. See Figure~\ref{fig:nreal}, the grey areas belong to the respective polytopes in case $\Theta^n_M$ is real.
\proofend

One may wonder whether Theorem~\ref{thm:chek} reflects a symplectic phenomenon or a smooth one. This is not obvious in general, but we discuss the case in which the moment polytope of $M$ is centrally symmetric. See Section~\ref{sec:intro} for a discussion and the classification of manifolds having this property. Although $\Theta^n_M$ is not real ($M$ has property $FS$ whenever $\Delta$ is centrally symmetric), we prove that it can be realized as the fixed point set of a smooth involution.

\begin{proposition}
\label{prop:smoothinv}
Let $M$ be a toric monotone symplectic manifold which has a centrally symmetric moment polytope $\Delta = -\Delta$. Then the Chekanov torus $\Theta^n_M$ is the fixed point set of a smooth involution.
\end{proposition}
\noindent

\proof It is proved in~\cite{BreKimMoo19} that the central fibre~$T_0$ is real whenever $\Delta = -\Delta$. Hence we can take an anti-symplectic involution $\sigma$ of $M$ such that $\Fix \sigma = T_0$. We claim that there is a $\psi \in \Diff(M)$ such that $\psi(\Theta^n_M) = T_0$. Then
\begin{equation*}
	\Theta^n_M = \Fix(\psi^{-1} \circ \sigma \circ \psi)
\end{equation*}
is the fixed point set of a smooth involution. The existence of $\psi$ follows from the proof of Lemma~\ref{lem:thetaversal}. Indeed, $\Theta^n_M$ is smoothly isotopic to all versal deformations $\Theta_{\mathbf{t},s}$ and whenever $\mathbf{t} \neq 0$, we have proved that $\Theta_{\mathbf{t},s}$ is isotopic to a toric fibre $T_x$. Since all toric fibres are isotopic, so are~$\Theta^n_M$ and~$T_0$.
\proofend

\subsection{More examples in $\times_n S^2$} In order to obtain more than only one example of non-real exotic Lagrangian torus in a given toric manifold, one may try to embed higher twist tori or products of Chekanov tori. We will discuss the second case here. For $\mathbf{k} = (k_1,\ldots,k_s)$ with $k_i \geqslant 2$ and $s \geqslant 1$, define the product
\begin{equation*}
	\label{eq:chekprod}
	\Theta^{\mathbf{k},m} = \Theta^{k_1} \times \ldots \times \Theta^{k_s} \times T_0^m \subset \CC^n, \quad 
 \sum_{i=1}^s k_i + m = n,
\end{equation*}
where~$T^m_0$ denotes the Clifford torus in~$\CC^m$. The image of such products under the standard moment map~$\nu$ in~$\CC^n$ is given by a hypercube formed by the product of diagonal segments
\begin{equation*}
	\nu(\Theta^{\mathbf{k},m}) = \{ (\underbrace{r_1,\ldots,r_1}_{k_1},\ldots,\underbrace{r_s,\ldots,r_s}_{k_s},\underbrace{0,\ldots,0}_m) \in \RR^n \vert \, \varepsilon -1 < r_i < \varepsilon \}.
\end{equation*}
In order to embed~$\Theta^{\mathbf{k},m}$ in a toric monotone symplectic manifold~$M$ with moment polytope~$\Delta$, one may try to apply the same strategy as for~$\Theta^n$, namely put~$\Delta$ in the normal form~\eqref{eq:ellnormal} and see if~$\nu(\Theta^{\mathbf{k},m})$ lies inside~$\Delta$. If it is so, the resulting torus is not real. 

\begin{proposition}
\label{prop:thetakmnr}
Let $M$ be a toric monotone symplectic manifold satisfying property $FS$. Assume furthermore that $\Theta^{\mathbf{k},m}$ can be embedded as described above. Then the image $\Theta^{\mathbf{k},m}_M \subset M$ is a monotone Lagrangian torus which is not real. 
\end{proposition}

\proof
Monotonicity follows from the same arguments as in the proof of Proposition~\ref{prop:Chekembed}. In order to prove that $\Theta^{\mathbf{k},m}_M$ is not real, we compute its displacement energy germ
\begin{equation}
\label{eq:germprod}
S_{\Theta_M^{\mathbf{k},m}} \stackrel{\otimes}{=} S_{T_0} \circ \phi_{\mathbf{k},m}.
\end{equation}
Here $\phi_{\mathbf{k},m} : \RR^n \rightarrow \RR^n$ is the piece-wise linear homeomorphism given as a product of the map $\phi$ defined as in the proof of Lemma~\ref{lem:thetaversal},
\begin{equation*}
\phi_{\mathbf{k},m} = \phi_{k_1} \times \ldots \times \phi_{k_s} \times \text{id}_m.
\end{equation*}
Indeed, note that the normal form~\eqref{eq:ellnormal} of~$v$ splits in~$\RR^n = \RR^{k_1} \times \ldots \times \RR^{k_s} \times \RR^m$ as the product of vertices in normal form. Hence the argument given in Lemma~\ref{lem:thetaversal} can be carried out on the factors. Let~$\pi: \RR^n \rightarrow \RR^{k_1}$ be the projection to the first~$k_1$ coordinates. Assume that the polytope~$\phi_{\mathbf{k},m}^{-1}(\Delta)$ is centrally symmetric. Then so is its projection~$\pi (\phi_{\mathbf{k},m}^{-1}(\Delta))$. By the product structure of~$\phi_{\mathbf{k},m}$, we have~$\pi (\phi_{\mathbf{k},m}^{-1}(\Delta)) = \phi_{k_1}^{-1}(\pi(\Delta))$. By convexity of~$\Delta$, the projection~$\pi(\Delta)$ is in normal form at the vertex~$\pi(v)$ and hence we can apply the same argument as in the proof of Theorem~\ref{thm:chek} to get a contradiction to the convexity of~$\pi(\Delta)$.
\proofend

To enumerate the non-real product tori $\Theta_M^{\mathbf{k}, m}$ that one obtains by this method up to symplectomorphism, one should now solve the following two problems: First, for which vertices $v$ does $\Theta^{\mathbf{k},m}$ fit into the normal form of $\Delta$ at $v$? Second, which of the so-obtained tori $\Theta^{\mathbf{k}, m}_M$ (also depending on the vertex~$v$) are exotic and which are pairwise non-symplectomorphic? For a general~$M$, both problems seem to involve complicated combinatorics outside of the scope of the present paper, whence we will only carry out the details for~$M=\times_n S^2$. In that case, all tori $\Theta^{\mathbf{k}, m}$ embed, the embedding does not depend on the vertex $v$, and all tori $\Theta^{\mathbf{k},m}_M$ turn out to be pairwise distinct.\\

Let $M = \times_n S^n$. As we have seen in~\eqref{eq:germlevel}, we can understand the versal deformation of $\Theta^n_M$ by applying~$\phi^{-1}$ to the moment polytope $\Delta = [-1,1]^n$ of $M$. We call the resulting polytopes \textbf{Chekanov polytopes} and denote them by 
\begin{equation*}
	\text{CP}_n = \phi^{-1}(\Delta).
\end{equation*}
We have a closer look at the geometry and the combinatorics of~$\text{CP}_n$. Notice that $s=\min\{x_1, \ldots,x_n\}$ is equal to~$-1$ on all facets that contain the vertex~$(-1,\ldots,-1)$. In other words, all of these facets are mapped to the hyperplane~$\{s=-1\}$ by $\phi^{-1}$. The one remaining vertex~$(1,\ldots,1)$ is mapped to~$e_n$. Hence $\text{CP}_n$ has the structure of a convex cone over the $(n-1)$-dimensional polytope~$P_{-1} = \{s=-1\} \cap \text{CP}_n$. In order to understand~$\text{CP}_n$, we thus need to understand $P_{-1}$ in the hyperplane $\{s=-1\} \cong \RR^{n-1}$. We claim that $P_{-1}$ is equal to the polytope obtaind by sweeping out the standard $(n-1)$-hypercube along~$r(-1,\ldots,-1)$ for all~$r \in [-1,1]$. This follows from equation~\eqref{eq:xt3}, which yields
\begin{equation*}
	P_{-1} = \{(x_1 - x_n, \ldots, x_{n-1} - x_n ) \, \vert \, x_i \in [-1,1] \text{ and } \min\{x_i\} = -1 \}.
\end{equation*}
The polytope $\text{CP}_n$ has $2^n - 1$ vertices, since $\phi^{-1}$ maps vertices to vertices except for $(-1,\ldots,-1)$ which is mapped to the interior of $P_{-1}$. The valencies of the vertices are given by 
\begin{equation}
	\label{eq:valencyvect}
	V(\text{CP}_n) = ((2^n-2)^{\times 1}, (n+1)^{\times (2^n - 2n -2)}, n^{ \times 2n}),
\end{equation}
for $n \geqslant 3$ where $l^{\times k}$ means that there are $k$ vertices with valency~$l$. For~$n=2$, we have $V(\text{CP}_2)=(2^{\times 3})$, as illustrated by Figure~\ref{fig:level3}. The general case can be seen as follows. The valency of the vertex at the apex of the cone is equal to the number of vertices of~$P_{-1}$ and hence equal to~$2^n -2$. We obtain the valency of any other vertex by adding $1$ to its valency when considered with respect to $P_{-1}$. Now let~$\epsilon = (\epsilon_1,\ldots,\epsilon_{n-1})$ for~$\epsilon_i=\pm 1$ be a vertex of the standard hypercube in $\RR^{n-1}$. Since $P_{-1}$ has central symmetry, we will restrict the count of valencies to the positive quadrant. Recall that $P_{-1}$ is obtained by sweeping out a copy of this hypercube centered in $(1,\ldots,1)$ along the vector~$(-1,\ldots,-1)$. Hence we can check what happens to the vertices of the standard hypercube under the sweeping along $r(-1,\ldots,-1)$ for small $r>0$. The vertex with all~$\epsilon_i = 1$ is untouched by this process and keeps valency~$n-1$. The vertex with all~$\epsilon_i = -1$ is erased by the sweeping. At all other vertices $\epsilon$, a new emanating edge is created by the sweeping, since whenever there is an~$\epsilon_i = -1$, the vector $\epsilon + t(-1,\ldots,-1)$ lies outside the hypercube for small $t>0$. Lastly, we check that exactly the~$n-1$ edges emanating from the vertex~$(-1,\ldots,-1)$ are deleted by the sweeping and hence there are~$n-1$ vertices whose valency decreases by one. Before sweeping, the facets containing the vertex $\epsilon$ have normal vectors
\begin{equation*}
	(\epsilon_1, 0, \ldots,0) , (0, \epsilon_2, \ldots, 0), \ldots, (0,0,\ldots,\epsilon_{n-1})
\end{equation*}
and a facet remains a facet after sweeping only if the projection of its normal vector onto~$(-1,\ldots,-1)$ is negative, i.e. if $\epsilon_i = 1$. Hence, in order for a given edge to be deleted by the sweeping, it has to be the intersection of $n-2$ facets with normal vectors defined by $\epsilon_i = -1$. These are precisely the edges emanating from $(-1,\ldots,-1)$. This yields the count in~\eqref{eq:valencyvect}.

\begin{proposition}
\label{prop:thetakm}
If two Lagrangian tori $\Theta^{\mathbf{k},m}_M$ and $\Theta^{\mathbf{k}',m'}_M$ in~$M=\times_n S^2$ are symplectomorphic, then $\mathbf{k} = \mathbf{k}'$ and $m = m'$.
\end{proposition}
\noindent 
Together with Proposition~\ref{prop:thetakmnr} we conclude that all the tori $\Theta_M^{\mathbf{k},m}$ in $\times_n S^2$ are not real and mutually not symplectomorphic. The number of such tori is $p(n)-1$, where $p(n)$ is the number of partitions of $n$.\\

\proofof{Proposition~\ref{prop:thetakm}}
By the product structure of $M$ and~\eqref{eq:germprod}, the level sets of the displacement energy germ of $\Theta^{\mathbf{k},m}_M$ are given by the product of Chekanov polytopes and intervals $I = [-1,1]$,
\begin{equation*}
	S_{\Theta^{\mathbf{k},m}_M}^{-1}(c) \stackrel{\otimes}{=} \text{CP}_{k_1} \times \ldots \times \text{CP}_{k_s} \times I^m.
\end{equation*}
Hence, it suffices to show that $\mathbf{k}=(k_1,\ldots,k_s)$ and $m$ are determined by the $\GL(n,\ZZ)$-equivalence class of $\times_i \text{CP}_{k_i} \times I^m$. In order to prove this, we associate to the latter polytopes the vector counting emanating edges at its vertices in decreasing order as in \eqref{eq:valencyvect}. This datum is a $\GL(n,\ZZ)$-invariant of polytopes. Note that if $P$ and $P'$ are polytopes, we have for the respective valency vectors
\begin{equation*}
	V(P \times P') = V(P) \oplus V(P'),
\end{equation*}
where the operation $\oplus$ on vectors $a = (a_1,\ldots,a_{k_1})$ and $b = (b_1,\ldots,b_{k_2})$ with $a_1 \geqslant a_2 \geqslant \ldots \geqslant a_{k_1}$ and $b_1 \geqslant b_2 \geqslant \ldots \geqslant b_{k_2}$ is defined as the vector of all possible sums in decrasing order
\begin{equation*}
	a \oplus b = (a_1 + b_1, \ldots, a_{k_1} + b_{k_2}).
\end{equation*}
This operation is commutative and associative, and hence we obtain
\begin{equation*}
	V(\text{CP}_{k_1} \times \ldots \times \text{CP}_{k_s} \times I^m) = V(\text{CP}_{k_1}) \oplus \dots \oplus V(\text{CP}_{k_s}) \oplus V(I^m).
\end{equation*}
Furthermore, this operation is invertible in the following sense. Let~$c = (c_1,\ldots,c_{k_1k_2})$ denote $a \oplus b$. Then $a$ is determined by $c$ and $b$; in other words, there is an operation $\ominus$ with $c \ominus b = a$. We will prove this by induction on the length $k_1$ of $a$. The case $k_1=1$ is obvious. In case $k_1=l+1$, note that $a_1$, $b_1$ and $c_1 = a_1 + b_1$ are by convention the maximal components of the corresponding vectors and hence $a_1$ is given by~$c_1 - b_1$. The situation can be reduced to the case~$k_1 = l$ by removing the value~$a_1$ from~$a$ and the values~$a_1 + b_1, \ldots, a_1 + b_{k_2}$ from~$c$.

We will now successively split off factors from the product polytope using the operation $\ominus$. First, notice that the multiplicity of the maximal entry of $V(\times_i \text{CP}_{k_i} \times I^m)$ determines $m$ and $p$, where $p$ is the number of times we have $k_i = 2$. Indeed, we have $V(I)=(1^{\times 2})$ and $V(\text{CP}_2)=(2^{\times 3})$ and by equation~\eqref{eq:valencyvect} the multiplicity of the maximal entry is given by~$2^m3^p$. Hence the prime decomposition of this multiplicity yields~$m$ and~$p$. After splitting off the corresponding factors, we can assume that~$m=0$ and~$k_i \geqslant 3$. Let $M_1$ and $M_2$ be the largest and the second largest component of the the valency vector. Then we have $M_1 = \sum_{i=1}^s 2^{k_i} -2s$ and $M_1 - M_2 = 2^{k_{\text{min}}} - k_{\text{min}} -3 $, where $k_{\text{min}}$ is minimal among all $k_i$. Therefore $M_1 - M_2$ determines $k_{\text{min}}$ and we can split off~$V(\text{CP}_{k_{\text{min}}})$ from the valency vector by using formula~\eqref{eq:valencyvect}.  
\proofend

\section{Appendix: Alternate approach using $J$-holomorphic disks}
\label{sec:appendix}

In this appendix, we outline an alternate approach to Theorem~\ref{t:toric} based on the count of~$J$-holomorphic Maslov~$2$ disks with boundary on the Lagrangian, which was introduced in~\cite{EliPol97} and~\cite{Che97} and was used in~\cite{Kim19b} to determine whether a given Lagrangian is real. This approach is less elementary than the above, but has the advantage of avoiding property~$FS$.\\

\noindent
Let~$T_0$ be the central fibre in a toric monotone symplectic manifold~$(M,\omega)$ with moment polytope~$\Delta = \{ \langle x , v_i \rangle \leqslant 1 \}$. Assuming that~$T_0$ is the fixed point set of an anti-symplectic involution~$\sigma$, we will show that~$\Delta$ is centrally symmetric. Fix an~$\omega$-compatible almost-complex structure~$J$ on~$M$ and a homology class~$\xi \in H_1(T_0,\ZZ)$ and define the moduli space
\begin{eqnarray*}
\mathcal{M}(T_0,J,\xi) = \{ u \colon (D, \partial D) \rightarrow (M,T_0) &\vert & u \; J\text{-holomorphic}, \\
&& \text{Maslov}(u) = 2, \\
&& [\partial u] = \xi \}/\sim,
\end{eqnarray*}
where~$\sim$ denotes the equivalence relation induced by reparametrizing the domain~$D$ by biholomorphisms fixing the point~$1 \in \partial D$. We can count (mod~$2$) the elements of $\mathcal{M}(T_0,J,\xi)$ whose boundary passes through a given point on~$T_0$ by taking the degree $n(T_0,J,\xi) \in \ZZ$ of the evaluation map
\begin{equation*}
\text{ev} \colon \mathcal{M}(T_0,J,\xi) \rightarrow T_0, \quad [u] \mapsto u(1).
\end{equation*}
See for example~\cite{Aur15} for details. As in~\cite{Kim19b}, we now assume in addition that~$\sigma^*J = -J$ and associate to every element~$[u] \in \mathcal{M}(T_0,J,\xi)$ its image under the anti-symplectic involution
\begin{equation*}
\mathcal{R} \colon [u] \mapsto [\sigma \circ u \circ \rho],
\end{equation*}
where~$\rho$ denotes complex conjugation on the disk. Note that the involution~$\mathcal{R}$ maps the moduli space~$\mathcal{M}(T_0,J,\xi)$ to~$\mathcal{M}(T_0,J,-\xi)$ since~$T_0$ is the fixed point set of~$\sigma$. By Cho and Oh~\cite{ChoOh06}, there exists a $J_0$-holomorphic disk in~$\mathcal{M}(T_0,J_0,\xi)$ if and only if~$\xi$ coincides with one of the primitive vectors~$v_i$ normal to the facets of the moment polytope~$\Delta$. Here,~$J_0$ denotes the K\"ahler complex structure. The regularity of~$J_0$ was shown in~\cite{ChoOh06} and that of~$J$ by Kim~\cite{Kim19b}. Hence the two counts~$n(T_0,J,\xi)$ and~$n(T_0,J_0,\xi)$ agree. Since the involution~$\mathcal{R}$ maps~$\mathcal{M}(T_0,J,\xi)$ to $\mathcal{M}(T_0,J,-\xi)$, we find that a given vector~$v$ appears as orthogonal vector to one of the facets if and only if~$-v$ does as well. Hence~$\Delta$ is invariant under central symmetry.

\begin{remark} \label{rk:wallcrossing}
{\rm
One can make a similar argument in the case of Chekanov tori by reformulating the information given by~$J$-holomorphic disks in terms of the Landau-Ginzburg potential (see~\cite{Aur09} or~\cite{PasTon20}). The so-called wall-crossing formulae describe how this potential behaves when passing from the Clifford to the Chekanov torus.
}
\end{remark}

\bibliographystyle{abbrv}
\bibliography{mybibfile}

\end{document}